\documentclass[a4paper,10pt]{amsart}

\usepackage[colorlinks=true,linkcolor=red,citecolor=blue,urlcolor=blue,hypertexnames=false]{hyperref}
\usepackage{bookmark}
\usepackage{amsthm,thmtools,amssymb,amsmath,amscd}

\usepackage{fancyhdr}
\usepackage{esint}
\usepackage{enumerate}

\usepackage{pictexwd,dcpic}

\swapnumbers
\declaretheorem[name=Theorem,numberwithin=section]{thm}
\declaretheorem[name=Remark,style=remark,sibling=thm]{rem}
\declaretheorem[name=Lemma,sibling=thm]{lemma}
\declaretheorem[name=Proposition,sibling=thm]{prop}

\declaretheorem[name=Assumption,style=definition,sibling=thm]{assum}

\numberwithin{equation}{section}

\usepackage{cleveref}
\crefname{lemma}{Lemma}{Lemmata}
\crefname{prop}{Proposition}{Propositions}
\crefname{thm}{Theorem}{Theorems}
\crefname{cor}{Corollary}{Corollaries}
\crefname{defn}{Definition}{Definitions}
\crefname{example}{Example}{Examples}
\crefname{rem}{Remark}{Remarks}
\crefname{assum}{Assumption}{Assumptions}
\crefname{notation}{Notation}{Notation}

\usepackage{autonum}

\newcommand{\ti}{\tilde}
\newcommand{\wt}{\widetilde}

\newcommand{\cn}{\colon}
\newcommand{\sub}{\subset}
\newcommand{\ov}{\overline}

\newcommand{\R}{\mathbb{R}}

\newcommand{\bbS}{\mathbb{S}}
\newcommand{\bbH}{\mathbb{H}}

\newcommand{\8}{\infty}

\newcommand{\al}{\alpha}
\newcommand{\be}{\beta}
\newcommand{\g}{\gamma}
\newcommand{\de}{\delta}
\newcommand{\e}{\epsilon}
\newcommand{\ka}{\kappa}

\newcommand{\s}{\sigma}

\newcommand{\p}{\varphi}

\newcommand{\vt}{\vartheta}

\newcommand{\G}{\Gamma}

\newcommand{\La}{\Lambda}

\newcommand{\cL}{\mathcal{L}}

\newcommand{\cW}{\mathcal{W}}

\newcommand{\cS}{\mathcal{S}}

\newcommand{\cF}{\mathcal{F}}

\newcommand{\del}{\partial}
\newcommand{\n}{\nabla}

\newcommand{\fa}{\forall}

\newcommand{\ip}[2]{\left\langle #1,#2 \right\rangle}
\newcommand{\fr}[2]{\frac{#1}{#2}}
\newcommand{\x}{\times}

\DeclareMathOperator{\id}{id}

\DeclareMathOperator{\tr}{tr}
\DeclareMathOperator{\Rm}{Rm}
\DeclareMathOperator{\Rc}{Rc}

\DeclareMathOperator{\grad}{grad}

\newcommand{\pf}[1]{\begin{proof} #1 \end{proof}}
\newcommand{\eq}[1]{\begin{equation}\begin{alignedat}{2} #1 \end{alignedat}\end{equation}}

\newcommand{\br}[1]{\left(#1\right)}

\newcommand{\Ra}{\Rightarrow}
\newcommand{\ra}{\rightarrow}
\newcommand{\hra}{\hookrightarrow}

\newcommand{\mrm}{\mathrm}

\newcommand{\hp}{\hphantom}
\newcommand{\q}{\quad}

\usepackage{tabularx}
\usepackage{multirow}
\usepackage{caption} 
\usepackage{makecell}
\usepackage{ mathrsfs }

\begin{document}

\title{Inverse curvature flows in Riemannian warped products}
\author{Julian Scheuer}
\subjclass[2010]{35J60, 53C21, 53C44, 58J05}
\keywords{Curvature flows, Inverse curvature flows, Warped products}
\date{\today}
\address{Julian Scheuer, Albert-Ludwigs-Universit{\"a}t, Mathematisches Institut, Eckerstr.~1,
79104 Freiburg, Germany}
\email{julian.scheuer@math.uni-freiburg.de}

\begin{abstract}
The long-time existence and umbilicity estimates for compact, graphical solutions to expanding curvature flows are deduced in Riemannian warped products of a real interval with a compact fibre. Notably we do not assume the ambient manifold to be rotationally symmetric, nor the radial curvature to converge, nor a lower bound on the ambient sectional curvature. The inverse speeds are given by powers $p\leq 1$ of a curvature function satisfying few common properties.
\end{abstract}
\maketitle

\section{Introduction}
This paper deals with expanding curvature flows of the form
\eq{\label{ICF}\dot{x}=\fr{1}{F^{p}}\nu,\q 0<p\leq1,}
where 
\eq{x\cn [0,T^{*})\x M^{n}\ra N^{n+1},\q n\geq 2,}
is a family of embeddings of a smooth, orientable, compact manifold $M^{n}$ and $N=N^{n+1}$ is a product
\eq{N=(R_{0},\8)\x \cS_{0}}
with metric
\eq{\bar{g}=dr^{2}+\vt^{2}(r)\s.}
Here $\vt\in C^{\8}((R_{0},\8))$ satisfies $\vt'>0$, $\vt''\geq 0$ and $(\cS_{0},\s)$ is a compact Riemannian manifold. In \eqref{ICF}, $F$ is a function evaluated at the Weingarten operator $\cW$ of the flow hypersurfaces $M_{t}=x(t,M)$ at the respective point $x$ and $\nu$ is the outward pointing normal, i.e.
\eq{\bar{g}(\nu,\del_{r})>0.} The detailed assumptions on the curvature function $F$ and on $N$ are the following.
 
\begin{assum}\label{ICF-F}
Let $\G\sub\R^{n}$ be an open, symmetric and convex cone containing the positive cone
\eq{\G_{+}=\{(\ka_{i})\in \R^{n}\cn \ka_{i}>0\q\fa 1\leq i\leq n \}} and suppose $f\in C^{\8}(\G)$ is a positive, symmetric, strictly monotone, $1$-homogeneous and concave function with
\eq{f(1,\dots,1)=n,\q f_{|\del\G}=0}
and associated curvature function $F=F(\cW)$, cf. \cref{CF}.
\end{assum}

Particular examples of curvature functions satisfying these assumptions are roots or quotients of elementary symmetric polynomials,
\eq{F=nH_{k}^{\fr 1k},\q F=n\fr{H_{k+1}}{H_{k}} }
and many more, cf. \cite{Andrews:/2007}.

In order to obtain good asymptotics we will make the following assumption on the warping function. This assumption will not be needed for the long-time existence.

\begin{assum}\label{ICF-N}
Assume the warping function $\vt\in C^{\8}((R_{0},\8))$ to satisfy
\eq{\limsup_{r\ra \8}\fr{\vt''\vt}{\vt'^{2}}<\8\q\text{and}\q \limsup_{\stackrel{r\ra \8}{\vt''(r)>0}}\fr{\vt'''\vt}{\vt'\vt''}<\8.}
\end{assum}

In the following theorem $\underline{\widehat{\Rc}}$ denotes the smallest eigenvalue of the Ricci tensor of $\s$ and $H_{k}$ denotes the curvature function determined by the $k$-th normalized elementary symmetric polynomial of the principal curvatures, compare \cref{CF} for further information. In this paper we aim to prove the following theorem.

\begin{thm}\label{ICF-Main}
Let $(\cS_{0},\s)$ be a smooth, compact and orientable Riemannian manifold of dimension $n\geq 2$, $R_{0}>0$, $N=(R_{0},\8)\x\cS_{0}$ and define a warped product metric on $N$,
\eq{\bar g=dr^{2}+\vt^{2}(r)\s,}
with $\vt\in C^{\8}((R_{0},\8))$,
 $\vt''\geq 0$ and $\vt'>0$. Let $0<p\leq 1$ and $F$ satisfy \cref{ICF-F}. Let 
\eq{x_{0}\cn M\hra N}
be the embedding of a hypersurface $M_{0}$, which is graphical over $\cS_{0}$, i.e. there exists $u\in C^{\8}(\cS_{0},(R_{0},\8))$ such that
\eq{M_{0}=\{(u(y),y)\cn y\in \cS_{0}\},} and such that all its $n$-tuples of principal curvatures belong to $\G$. 
\begin{enumerate}[(i)] \item Assume either of the following properties to hold:
\eq{
	&(a)~\s~ \text{has non-negative sectional curvature.}\\
	&(b)~F=n\fr{H_{k+1}}{H_{k}},\q 0\leq k\leq n-1. }
Then there exists a unique immortal solution 
\eq{x\cn [0,\8)\x M&\ra N}
of
\eq{\label{ICF-Ini}\dot{x}&=\fr{1}{F^{p}}\nu\\
		x(0,\cdot)&=x_{0},}
which is also graphical over $\cS_{0}$, i.e. $\ip{\nu}{\del_{r}}>0$. 
\item Assume $\s$ has non-negative sectional curvature and each of the following properties:
\eq{
 &(A)~\cref{ICF-N}~\text{holds}.\\
&(B)~\sup_{r>0}\vt'(r)<\8~\text{and}~p=1\q\Ra\q \underline{\widehat{\Rc}}>0~\text{and}~  F=n\fr{H_{k+1}}{H_{k}}, \q 0\leq k\leq n-1.\\
&(C)~\sup_{r>0}\vt'(r)=\8~\text{and}~p=1\q\Ra\q \liminf_{r\ra \8}\fr{\vt''\vt}{\vt'^{2}}>0.
}
Then the flow hypersurfaces become umbilical with the rate
\eq{\label{ICF-Umb}\left|h^{i}_{j}-\fr{\vt'}{\vt}\de^{i}_{j}\right|\leq ct\fr{\vt'^{1-p(p+1)}}{\vt},}
where the $t$-factor may be dropped in case $p<1$ or bounded $\vt'$ and may even be replaced by $e^{-\al t}$ for some positive $\al$ if $\vt'$ is bounded and $p=1$.
\end{enumerate}
\end{thm}

Let us make some remarks on the technical assumptions made in \cref{ICF-Main}.

\begin{rem}
\begin{enumerate}[(i)]
\item The assumptions in statement (i) of \cref{ICF-Main} are optimal in the sense, that for example in a spherical ambient space with $\vt''<0$ the inverse mean curvature flow only exists for a finite time, cf. \cite{Gerhardt:/2015,MakowskiScheuer:11/2016} and for $p>1$ the maximal existence is finite if $N=\R^{n+1}$, cf. \cite{Gerhardt:01/2014}.
\item The assumption on the sectional curvature of $\s$ can be relaxed. The crucial point, where we use this assumption is in the first gradient estimates, especially in estimate \eqref{ICF-grad-bound-2}, where we throw away the term involving $\widehat{\Rm}$, if $F$ is general. However, under a further suitable technical assumption we could also absorb it into the first line of this equation. For the special case of the inverse mean curvature flow in the Reissner-Nordstr\"om manifolds this has been accomplished in the recent preprint \cite{ChenLiZhou:10/2017}. However, in order to avoid too many technical assumptions, we will not improve the main result in this direction here, except that we prove the long-time existence in general, provided that $F$ is a quotient of the $H_{k}$. For the IMCF this was also accomplished in \cite{LiWei:/2017,Zhou:06/2017}.
\item The rates of convergence in this theorem can be improved, if the ambient sectional curvatures approach each other at infinity. Such results have been accomplished for example in \cite{ChenMao:04/2017,Lu:09/2016,Scheuer:01/2017} in case $p=1$ and in \cite{Scheuer:05/2015} in case $p<1$ in the hyperbolic space. Since the main aim of this work is to deal with spaces in which the limits of the quantities in \cref{ICF-N} do not exist (if $\s$ is the round metric this implies that $N$ is not asymptotically a spaceform), we will not pursue these optimal estimates here and stick to the best we could accomplish in general ambient spaces. To the best of my knowledge, the only result in such general spaces is the analogous result for the inverse mean curvature flow proven in \cite{Scheuer:01/2017}.
\item The question, whether \eqref{ICF-Umb} implies that the flow hypersurfaces do become almost umbilical, depends on the ambient space $N$ and on $p$. However, if $p=1$, the analysis in \cite[Prop.~3.1]{Scheuer:01/2017} implies that $\vt$ grows exponentially. Hence in this case we obtain exponential decay of $\cW-\fr{\vt'}{\vt}\id$.
\item In case $p=1$, the gradient decay estimates obtained in \cref{ICF-optGrad} are optimal even if the ambient space is asymptotically a spaceform. Compare the explanation in \cite[Rem.~1.5]{Scheuer:01/2017}.
\item In case $p=1$ the estimate \eqref{ICF-Umb} turned out to be strong enough to obtain geometric inequalities, for example in \cite{BrendleHungWang:01/2016,GeWangWuXia:03/2015,LuMiao:03/2017,Wang:/2015}. We are optimistic that \cref{ICF-Main} will be helpful with such applications as well.
\end{enumerate}
\end{rem}

The motivation to analyse the behaviour of inverse curvature flows has mostly been driven by their power to deduce geometric inequalities for hypersurfaces. The most prominent example is the proof of the Riemannian Penrose inequality due to Huisken/Ilmanen \cite{HuiskenIlmanen:/2001}, building on the observation made by Geroch \cite{Geroch:12/1973} and Jang/Wald \cite{JangWald:01/1977} that the Hawking mass of a connected surface is non-decreasing under the inverse mean curvature flow (IMCF) with $F=H$ and $p=1$, if the ambient scalar curvature is non-negative.
Since for general initial data the IMCF may develop singularities, Huisken and Ilmanen defined a notion of a weak solution for this flow, maintaining the Geroch monotonicity. This enabled them to prove the Riemannian Penrose inequality. For a short outline of their procedure also compare \cite{HuiskenIlmanen:/1997}.

Also the classical solution to IMCF has lead to very interesting applications. A crucial feature of this flow in $\R^{n+1}$ is that one does not need to require convexity of the initial hypersurface to avoid finite time singularities. Namely, Gerhardt \cite{Gerhardt:/1990} and Urbas \cite{Urbas:/1990} proved the long-time existence even for more general flows in $\R^{n+1}$,
with $F$ satisfying \cref{ICF-F}, $p=1$ and a starshaped initial hypersurface $M_{0}$ with $F_{|M_{0}}>0$. Furthermore, after exponential rescaling, the flow converges to a sphere smoothly. This result, with $F=nH_{k+1}/H_{k}$, was later exploited by Guan/Li \cite{GuanLi:08/2009} to generalise the Alex\-androv-Fenchel quermassintegral inequalities
 from the convex setting to the starshaped and $H_{k+1}$-convex setting. Since then a cascade of similar results followed by the same method (monotone quantity plus some convergence result) in various ambient spaces. The tough parts are to find the monotone quantity and to prove a sufficient convergence result. Examples of other results in this direction are a generalised Minkowski-type inequality in the anti-de Sitter-Schwarzschild manifold due to Brendle/Hung/Wang \cite{BrendleHungWang:01/2016}, Alexandrov-Fenchel-type inequalities in the hyperbolic space \cite{De-LimaGirao:04/2016,GeWangWu:04/2014,Hu:10/2016,LiWeiXiong:03/2014,WeiXiong:/2015} and in the sphere \cite{GiraoPinheiro:12/2017,MakowskiScheuer:11/2016,WeiXiong:/2015}. Further similar applications can be found in \cite{GeWangWuXia:03/2015,KwongMiao:/2014,LuMiao:03/2017,Wang:/2015}.
 
In many of these papers, there was a need to investigate the asymptotical behaviour of the corresponding inverse curvature flow separately, since a unified treatment had not been present. Hence, a branch of research solely dealing with inverse curvature flows has developed within the community, where the main aims are to generalise the convergence results in various directions (concerning flow speed and ambient space). A step towards generalising the ambient space was made by the author with the paper \cite{Scheuer:01/2017}, where the IMCF was considered in rotationally symmetric warped products under assumptions similar to \cref{ICF-N}. Before (and after) that, some more special ambient spaces were treated, which, to the best of my knowledge, all assumed convergence of the quantities in \cref{ICF-N}. Instead of giving a description of the available results verbally, the following table is supposed to give an overview as broad as I could accomplish over the previous results on {\it{smooth, inverse curvature flows of closed hypersurfaces in Riemannian warped products}}. The topics they cover are for example long-time existence, asymptotic behaviour, solitons and others. We point out that, in order to keep things manageable, we leave aside treatments of contracting flows, weak solutions, flows in Lorentzian manifolds, flows of entire graphs, flows with boundary conditions, anisotropic flows and flows with constraints (e.g. volume preserving flows). 

\renewcommand{\arraystretch}{1.3}
\vspace{10pt}
    \begin{tabular}{|*{5}{c|}}
    \hline
      N/F &   & $F=n\fr{H_{k+1}}{H_{k}}$ & \makecell{$F$ more general\\ and $p=1$} & \makecell{$p\neq 1$ or\\ non-hom. speed} \\ \hline
\multirow{3}{*}{ CSC}   & $\R^{n+1}$  &\makecell{\cite{CastroLerma:06/2017,ChowLiouTsai:09/1996,DruganLeeWheeler:/2016}\\\cite{HuiskenIlmanen:/2008}} &  \makecell{\cite{ChowChowFong:01/2017,Gerhardt:/1990,Li:06/2010}\\\cite{Smoczyk:/2005,Urbas:/1990,Urbas:/1991} } & \makecell{\cite{Andrews:09/1994,Andrews:/1998,BIS4,ChowTsai:/1996,ChowTsai:/1997}\\\cite{ChowTsai:/1998,Gerhardt:01/2014,Ivaki:03/2015,IvochkinaNehringTomi:/2000}\\\cite{KronerScheuer:03/2017,LiWangWei:09/2016,Li:/2011,LinTsai:02/2010}\\\cite{Scheuer:07/2016,Schnuerer:/2006,Wei:09/2017}} \\ \cline{2-5}
   & $\bbH^{n+1}$ & \cite{Ding:01/2011,HungWang:09/2015}    &\makecell{\cite{Gerhardt:11/2011,Liu:05/2017,Yu:/2017}} & \cite{LiWangWei:09/2016,Scheuer:01/2015,Scheuer:05/2015,Wei:09/2017}  \\ \cline{2-5}
   & $\bbS^{n+1}$ & & \cite{Gerhardt:/2015,Liu:05/2017} & \cite{BIS1,BIS4,MakowskiScheuer:11/2016,Wei:09/2017} \\ \hline
\multirow{2}{50pt}{\makecell{Asympt.\\ CSC}}   & $\R^{n+1}$ & \cite{Ding:01/2011,LiWei:/2017,LuMiao:03/2017}  & &  \\ \cline{2-5}
   & $\bbH^{n+1}$ & \cite{BrendleHungWang:01/2016,Lu:09/2016,Neves:/2010}   & \cite{ChenMao:04/2017} &  \\ \hline
   \multirow{2}{*}{\makecell{More\\ general}} & \makecell{$\fr{\vt''\vt}{\vt'^{2}}$\\ converges}  & \makecell{\cite{ChenLiZhou:10/2017,LiWei:/2017,Mullins:10/2016,Zhou:06/2017}\\ } & &  \\     \cline{2-5}
   	& \makecell{\ref{ICF-N}} & \cite{Scheuer:01/2017} & & \\ \hline
   
    \end{tabular}
  
    \vspace{10pt}
\noindent
\renewcommand{\arraystretch}{1}
Note that a reference only appears in the most general slot it can be placed. Also note that there are few works on the inverse mean curvature flow in ambient spaces which are not warped products, \cite{Allen:08/2017,KoikeSakai:/2015,Pipoli:10/2016,Pipoli:04/2017}.
This paper aims to fill some gaps in this table, especially in the two bottom rows, and is organised as follows. \Cref{Prelim} collects some notation, conventions, basic facts about curvature functions and the relevant evolution equations. In \cref{LTE} we treat the long-time existence and in \cref{Asym} we analyse the asymptotic behaviour and finish the proof of \cref{ICF-Main}.

\section{Preliminaries}\label{Prelim}

\subsection{Notation and conventions}\label{Conv}

In this paper we deal with embedded hypersurfaces
\eq{x\cn M\hra N}
of a smooth, closed and orientable manifold $M^n$ into an ambient Riemannian manifold $(N^{n+1},\bar{g})$. All geometric quantities of $N$ will be furnished with an overbar, e.g. $\bar g=(\bar{g}_{\al\be})$ for the metric, $\bar\n$ for its Levi-Civita connection etc. In coordinate expressions, greek indices run from $0$ to $n$. For the quantities induced by the embedding $x$, we use latin indices running from $1$ to $n$, e.g. for the induced metric $g=(g_{ij})$ with Levi-Civita connection $\n$.
For a $(k,l)$ tensor field $T$ on $M$, its covariant derivative $\n T$ is a $(k,l+1)$ tensor field given by
\eq{&(\n T)(Y^1,\dots, Y^k,X_1,\dots,X_l,X)\\
=~&(\n_{X} T)(Y^1,\dots, Y^k,X_1,\dots,X_{l})\\
						=~&X(T(Y^1,\dots,Y^k,X_1,\dots,X_l))-T(\n_X Y^1,Y^2,\dots, Y^k,X_1,\dots,X_l)-\ldots\\
                        \hp{=}~&-T(Y^1,\dots,Y^k,X_1,\dots,X_{l-1}\n_X X_l),}
the coordinate expression of which is denoted by
\eq{\n T=\br{T^{i_1\dots i_k}_{j_1\dots j_l;j_{l+1}}}.}
The index appearing after the semicolon indicates the derivative index.

Our convention for the $(1,3)$-Riemannian curvature tensor $\Rm$ of a connection $\n$ is
\eq{\Rm(X,Y)Z=\n_{X}\n_{Y}Z-\n_{Y}\n_{X}Z-\n_{[X,Y]}Z,} 
where $X,Y,Z$ are vector fields and where $[X,Y]$ is the Lie-bracket
\eq{[X,Y]\p=X(Y \p)-Y(X \p)\q\fa \p\in C^{\8}(M).}
The purely covariant Riemannian curvature tensor is defined by lowering to the fourth slot:
\eq{\Rm(X,Y,Z,W)=g(\Rm(X,Y)Z,W).}
Finally the Ricci curvature is
\eq{\Rc(X,Y)=\tr\br{\Rm(\cdot,X)Y}.}
For metrics $(g_{ij})$ we always denote its dual by $(g^{ij})$, i.e.
\eq{\de^{i}_{j}=g^{ik}g_{kj}.}

 The induced geometry of $M$ is governed by the following relations. The second fundamental form $h=(h_{ij})$ is defined by the Gaussian formula
\eq{\label{GF}\bar\n_{X}Y=\n_X Y-h(X,Y)\nu,}
where $\nu$ is a normal field. The Weingarten endomorphism $\cW=(h^i_j)$ is defined by $h^i_j=g^{ki}h_{kj}$ and we have the Weingarten equation
\eq{\label{Weingarten} \bar\n_X \nu=\cW(X). }
We also have the Codazzi equation
\eq{\label{Codazzi}\n_Z h(X,Y)-\n_Y h(X,Z)=-\ov{\Rm}(\nu,X,Y,Z). }
Let us record this equation is coordinates:
\eq{\label{Codazzi-1}h_{ij;k}-h_{ik;j}=-\ov{\Rm}(\nu,x_{;i},x_{;j},x_{;k}).}

The Gauss equation states
\eq{\label{GE}\Rm(W,X,Y,Z)=\ov{\Rm}(W,X,Y,Z)+h(W,Z)h(X,Y)-h(W,Y)h(X,Z)}
or in coordinates
\eq{R_{ijkl}=\ov{\Rm}(x_{;i},x_{;j},x_{;k},x_{;l})+h_{il}h_{jk}-h_{ik}h_{jl}.}

\subsection*{Warped products}
Throughout this paper we assume that the ambient manifold is a {\it{warped product}} of the form
\eq{\label{N}(N,\bar g)=(I\x \cS_{0},\bar g),}
where $I=(R_{0},\8)$, $(\cS_{0},\s)$ is an $n$-dimensional compact Riemannian manifold and
\eq{\label{Warped}\bar g=dr^{2}+\vt^{2}(r)\s}
with $\vt\in C^{\8}((R_{0},\8))$. We will need to know how the curvature tensor of $\bar g$ arises from the curvature tensors of $dr^{2}$ and $\s$. The relevant formulae can be found in \cite[Ch.~7, Prop.~42]{ONeill:/1983}. We state them here for further use, but adapted to our curvature convention, which differs from the one in op.~cit. We denote by $\mathscr{L}(\R)$ and $\mathscr{L}(\cS_{0})$ the space of all vector field on $\R$ resp. $\cS_{0}$ lifted to $N$.

\begin{lemma}\textsc{(\cite[Ch.~7,~Prop.~42]{ONeill:/1983})}
Let $N$ be given as above. If  $X,Y,Z\in \mathscr{L}(\R)$ and $U,V,W\in \mathscr{L}(\cS_{0})$, then the Riemannian curvature tensor of $N$ is given by
\begin{enumerate}
\item[(i)] $\ov{\Rm}(X,Y)Z=0,$
\item[(ii)] $\ov{\Rm}(V,X)Y=-\fr{\bar{\n}^{2}\vt(X,Y)}{\vt}V=-\fr{\vt''}{\vt}\bar{g}(X,Y)V,$
\item[(iii)] $\ov{\Rm}(X,Y)V=\ov{\Rm}(V,W)X=0,$
\item[(iv)] $\ov{\Rm}(X,V)W=-\fr{\vt''}{\vt}\bar g(V,W)X$
\item[(v)] $\ov{\Rm}(V,W)U=\widetilde{\Rm}(V,W)U-\fr{\vt'^{2}}{\vt^{2}}(\bar g(W,U)V-\bar{g}(V,U)W),$
\end{enumerate}
where $\widetilde{\Rm}$ is the lift of the Riemann tensor of the fibre $(\cS_{0},\vt^{2}(r)\s)$ under the projection $\pi\cn N\ra \cS_{0}$.
\end{lemma}

It will turn out to be convenient to have a closed coordinate expression for $\ov{\Rm}$, which follows easily from checking all of the five cases.

\begin{lemma}
In coordinates the Riemannian curvature tensor of the warped product
\eq{(N,\bar g)=(I\x \cS_{0},dr^{2}+\vt^{2}(r)\s)}
is given by
\eq{\label{Warped-Rm}{\bar{R}_{\al\be\g}}^{\e}&=\br{\br{\fr{\vt''}{\vt}-\fr{\vt'^{2}}{\vt^{2}}}
\bar{S}_{\al'\be'\g'\de'}+\wt{R}_{\al'\be'\g'\de'}}P^{\al'}_{\al}P^{\be'}_{\be}P^{\g'}_{\g}P^{\de'\e}-\fr{\vt''}{\vt}\bar{S}_{\al\be\g}^{\e},
}
where 
\eq{\bar{S}_{\al\be\g}^{\e}=\bar{g}_{\be\g}\de^{\e}_{\al}-\bar{g}_{\al\g}\de^{\e}_{\be}}
and
\eq{P=\id-\fr{\del}{\del r}\otimes dr.}
\end{lemma}

Hence we obtain a formula for the derivative of $\ov{\Rm}$.

\begin{lemma}\label{Warped-Rm-Der}
The coordinate functions of the covariant derivative of the $(0,4)$-curvature tensor are given by
\eq{\label{Warped-Rm-Der-1}\bar{R}_{\al\be\g\de;\e}&=-\br{\fr{\vt''}{\vt}}'r_{;\e}\bar{S}_{\al\be\g\de}+\br{\fr{\vt''}{\vt}-\fr{\vt'^{2}}{\vt^{2}}}'r_{;\e}\bar{S}_{\al'\be'\g'\de'}P^{\al'}_{\al}P^{\be'}_{\be}P^{\g'}_{\g}P^{\de'}_{\de}\\
	&\hp{=}+\wt{R}_{\al'\be'\g'\de';\e}P^{\al'}_{\al}P^{\be'}_{\be}P^{\g'}_{\g}P^{\de'}_{\de}-\fr{\vt'}{\vt}r_{;\al}\bar{T}_{\al'\be'\g'\de'}P^{\al'}_{\e}P^{\be'}_{\be}P^{\g'}_{\g}P^{\de'}_{\de}\\
	&\hp{=}-\fr{\vt'}{\vt}r_{;\be}\bar{T}_{\al'\be'\g'\de'}P^{\al'}_{\al}P^{\be'}_{\e}P^{\g'}_{\g}P^{\de'}_{\de}-\fr{\vt'}{\vt}r_{;\g}\bar{T}_{\al'\be'\g'\de'}P^{\al'}_{\al}P^{\be'}_{\be}P^{\g'}_{\e}P^{\de'}_{\de}\\
	&\hp{=}-\fr{\vt'}{\vt}r_{;\de}\bar{T}_{\al'\be'\g'\de'}P^{\al'}_{\al}P^{\be'}_{\be}P^{\g'}_{\g}P^{\de'}_{\e},}
where
\eq{\bar{T}_{\al'\be'\g'\de'}=\br{\fr{\vt''}{\vt}-\fr{\vt'^{2}}{\vt^{2}}}
\bar{S}_{\al'\be'\g'\de'}+\ti{R}_{\al'\be'\g'\de'}.}
	\end{lemma}

\pf{
Denote by $\bar\G^{\g}_{\al\be}$ the Christoffel symbols, i.e.
\eq{\bar\G^{\g}_{\al\be}\fr{\del}{\del x^{\g}}=\bar{\n}_{\fr{\del}{\del x^{\al}}}\fr{\del}{\del x^{\be}}}
and
\eq{\bar\G^{\g}_{\al\be}=\fr{1}{2}\bar{g}^{\g\de}\br{\fr{\del}{\del x^{\be}}\bar{g}_{\al\de}+\fr{\del}{\del x^{\al}}\bar{g}_{\be\de}-\fr{\del}{\del x^{\de}}\bar{g}_{\al\be}}.}
Using the definition of the metric we see
\eq{\bar{\G}^{0}_{\al\e}=-\fr{\vt'}{\vt}\bar{g}_{\al'\be'}P^{\al'}_{\al}P^{\be'}_{\e}=-\fr{\vt'}{\vt}\bar{g}_{\al'\e}P^{\al'}_{\al},\q \bar{\G}^{\al'}_{0\e}=\fr{\vt'}{\vt}P^{\al'}_{\e}}
and hence there holds
\eq{P^{\al'}_{\al;\e}=-r^{\al'}_{;\e}r_{\al}-{r_{;}}^{\al'}r_{;\al\e}&=-\bar{\G}^{\al'}_{0\e}r_{;\al}+{r_{;}}^{\al'}\bar{\G}^{0}_{\al\e}\\
			&=-\fr{\vt'}{\vt}P^{\al'}_{\e}r_{;\al}-\fr{\vt'}{\vt}{r_{;}}^{\al'}\bar{g}_{\g'\e}P^{\g'}_{\al}.}
There holds
\eq{\bar{T}_{\al'\be'\g'\de'}{r_{;}}^{\al'}P_{\e\al}P^{\be'}_{\be}P^{\g'}_{\g}P^{\de'}_{\de}=0}
and hence differentiation of \eqref{Warped-Rm} gives
\eq{\bar{R}_{\al\be\g\de;\e}&=-\br{\fr{\vt''}{\vt}}'r_{;\e}\bar{S}_{\al\be\g\de}+\br{\fr{\vt''}{\vt}-\fr{\vt'^{2}}{\vt^{2}}}'r_{;\e}\bar{S}_{\al'\be'\g'\de'}P^{\al'}_{\al}P^{\be'}_{\be}P^{\g'}_{\g}P^{\de'}_{\de}\\
	&\hp{=}+\ti{R}_{\al'\be'\g'\de';\e}P^{\al'}_{\al}P^{\be'}_{\be}P^{\g'}_{\g}P^{\de'}_{\de}-\fr{\vt'}{\vt}r_{;\al}\bar{T}_{\al'\be'\g'\de'}P^{\al'}_{\e}P^{\be'}_{\be}P^{\g'}_{\g}P^{\de'}_{\de}\\
	&\hp{=}-\fr{\vt'}{\vt}r_{;\be}\bar{T}_{\al'\be'\g'\de'}P^{\al'}_{\al}P^{\be'}_{\e}P^{\g'}_{\g}P^{\de'}_{\de}-\fr{\vt'}{\vt}r_{;\g}\bar{T}_{\al'\be'\g'\de'}P^{\al'}_{\al}P^{\be'}_{\be}P^{\g'}_{\e}P^{\de'}_{\de}\\
	&\hp{=}-\fr{\vt'}{\vt}r_{;\de}\bar{T}_{\al'\be'\g'\de'}P^{\al'}_{\al}P^{\be'}_{\be}P^{\g'}_{\g}P^{\de'}_{\e},}
which is the claimed formula.
}

 We will later have to deal with the $\ti{R}_{\al'\be'\g'\de';\e}$-term in \eqref{Warped-Rm-Der-1}.

\begin{lemma}\label{tildeRm}
For every $r_{0}>R_{0}$ there exists a constant $c$ such that
\eq{\|\bar\n\wt{\Rm}\|\leq c\fr{\vt'}{\vt^{3}}.}
\end{lemma}

\pf{We define a $\bar{g}$-orthonormal frame $(\ti{e}_{\al})_{0\leq \al\leq n}$ as follows: 
\eq{e_{0}=\ti{e}_{0}=\del_{r}}
and, given a $\s$-orthonormal frame $(e_{i})_{1\leq i\leq n}$ on $\cS_{0}$ we put
\eq{\ti{e}_{i}=\vt^{-1}e_{i}.}
Then clearly
\eq{\bar{g}(\ti{e}_{\al},\ti{e}_{\be})=\de_{\al\be},\q 0\leq \al,\be\leq n.}
To prove the lemma, it suffices to estimate the components of $\bar\n\wt{\Rm}$ accordingly with respect to this frame. There holds
\eq{\label{tildeRm-1}\bar{\n}_{\ti{e}_{\e}}\wt{\Rm}(\ti{e}_{\al},\ti{e}_{\be},\ti{e}_{\g},\ti{e}_{\de})&=\ti{e}_{\e}\br{\wt{\Rm}(\ti{e}_{\al},\ti{e}_{\be},\ti{e}_{\g},\ti{e}_{\de})}-\wt{\Rm}(\bar\n_{\ti{e}_{\e}}\ti{e}_{\al},\ti{e}_{\be},\ti{e}_{\g},\ti{e}_{\de})\\
		&\hp{=}-\wt{\Rm}(\ti{e}_{\al},\bar{\n}_{\ti{e}_{\e}}\ti{e}_{\be},\ti{e}_{\g},\ti{e}_{\de})-\wt{\Rm}(\ti{e}_{\al},\ti{e}_{\be},\bar{\n}_{\ti{e}_{\e}}\ti{e}_{\g},\ti{e}_{\de})\\
		&\hp{=}-\wt{\Rm}(\ti{e}_{\al},\ti{e}_{\be},\ti{e}_{\g},\bar{\n}_{\ti{e}_{\e}}\ti{e}_{\de}).}
There holds
\eq{\wt{\Rm}(\ti{e}_{\al},\ti{e}_{\be},\ti{e}_{\g},\ti{e}_{\de})=\vt^{-2}\widehat{\Rm}(\pi_{\ast}e_{\al},\pi_{\ast}e_{\be},\pi_{\ast}e_{\g},\pi_{\ast}e_{\de}),}
where $\widehat{\Rm}$ is the Riemann tensor of $\s$. Hence
\eq{\label{tildeRm-2}\ti{e}_{\e}\br{\wt{\Rm}(\ti{e}_{\al},\ti{e}_{\be},\ti{e}_{\g},\ti{e}_{\de})}=\begin{cases} -\fr{2\vt'}{\vt^{3}} \widehat{\Rm}(\pi_{\ast}e_{\al},\pi_{\ast}e_{\be},\pi_{\ast}e_{\g},\pi_{\ast}e_{\de}), &\e=0\\
			\vt^{-3}e_{\e}\br{\widehat{\Rm}(\pi_{\ast}e_{\al},\pi_{\ast}e_{\be},\pi_{\ast}e_{\g},\pi_{\ast}e_{\de})}, &\e\neq 0. \end{cases}}
From \cite[Ch.~7, Prop.~35]{ONeill:/1983} we obtain
\eq{\label{tildeRm-3}\pi_{\ast}\bar{\n}_{\ti{e}_{\e}}\ti{e}_{\al}=\begin{cases}\fr{\vt'}{\vt}\pi_{\ast}\ti{e}_{\al} ,&\e=0\\
\hat{\n}_{\ti{e}_{\e}}\ti{e}_{\al}, & \e\neq 0, \al\neq 0\\
	\fr{\vt'}{\vt}\ti{e}_{\e}, & \e\neq 0, \al=0, \end{cases}}
where $\hat{\n}$ is the Levi-Civita connection of $\s$. In case $\e\neq 0$, $\al\neq 0$ we have
\eq{\hat{\n}_{\ti{e}_{\e}}\ti{e}_{\al}=\vt^{-1}\hat{\n}_{e_{\e}}(\vt^{-1}e_{\al})=\vt^{-2}\hat{\n}_{e_{\e}}e_{\al}}
and
\eq{\wt{\Rm}(\bar\n_{\ti{e}_{\e}}\ti{e}_{\al},\ti{e}_{\be},\ti{e}_{\g},\ti{e}_{\de})=\vt^{-3}\widehat{\Rm }(\hat{\n}_{e_{\e}}e_{\al},e_{\be},e_{\g},e_{\de}).}
Using \eqref{tildeRm-2} and \eqref{tildeRm-3} in \eqref{tildeRm-1} in any of the cases, we obtain the desired estimate, since $\vt'\geq c_{r_{0}}>0$ on every interval $[r_{0},\8)$, giving the estimate
\eq{\vt^{-3}\leq c\fr{\vt'}{\vt^{3}}.}
}

\begin{rem}

For example, if $\s$ is the round metric on $\cS_{0}=\bbS^{n}$, then
\eq{\ti{R}_{\al\be\g\de}=\fr{1}{\vt^{2}}\bar{S}_{\al'\be'\g'\de'}P^{\al'}_{\al}P^{\be'}_{\be}P^{\g'}_{\g}P^{\de'}_{\de}.}
\end{rem}

\subsubsection*{Graphs in warped products}
 
 The hypersurfaces
\eq{x\cn M\hra N}
we deal with in this paper will all be graphs over $\cS_{0},$
\eq{x(M)=\{(u(y),y)\cn y\in \cS_{0}\}=\{(u(y(\xi)),y(\xi))\cn \xi \in M\},}
where 
\eq{u\cn \cS_{0}\ra (R_{0},\8)}
is smooth. Along $M$ we will always use the {\it{outward}} pointing normal 
\eq{\nu=v^{-1}(1,-\vt^{-2}\s^{ik}u_{;k}),}
where 
\eq{v^2=1+\vt^{-2}\s^{ij}u_{;i}u_{;j},}
and use this normal in the Gaussian formula \eqref{GF}. The support function of $M$ is defined by
\eq{\label{support}s=\bar{g}(\vt\del_r,\nu)=\fr{\vt}{v}.}
There is a relation between the second fundamental form and the graph function on the hypersurface. Let 
\eq{\bar h=\vt'\vt\s,}
then there holds
\eq{\label{graph-h2}v^{-1}h_{ij}=-u_{;ij}+\bar h_{ij},}
cf. 
\cite[equ.~(1.5.10)]{Gerhardt:/2006}. The induced metric is given by
\eq{g_{ij}=u_{;i}u_{;j}+\vt^2\s_{ij}}
and hence
\eq{\label{graph-h}v^{-1}h_{ij}=-u_{;ij}+\fr{\vt'}{\vt}g_{ij}-\fr{\vt'}{\vt}u_{;i}u_{;j}.}
In order to deduce the gradient estimates, it has proven to be useful to consider the function 
\eq{\p\cn \cS_{0}\ra \R }
 \eq{\label{varphi}\varphi(y)=\int_{\inf u_{0}}^{u(y)} \frac{1}{\vt(s)}~ds.} There holds
 \eq{\label{Wphi}h^{j}_{i}=\frac{\vt'}{\vt v}\de_i^j-\frac{1}{\vt v}\tilde{g}^{jk}\varphi_{:ki},}
where \eq{\tilde{g}^{ij}=\sigma^{ij}-\frac{\varphi_{:}^{\ i}\varphi_{:}^{\ j}}{v^2}} and the covariant derivative and index raising is performed with respect to $\s$, cf. \cite[equ.~(3.26)]{Gerhardt:11/2011}. We will use $\hat \n$ to denote the covariant derivative on $\cS_{0}$ throughout this paper.

\subsection{Curvature functions}\label{CF}
Let $\G\sub\R^{n}$ be an open and symmetric cone. In \cref{ICF-F} the symmetric function $f\in C^{\8}(\G)$ is supposed to be evaluated at the principal curvatures of the flow hypersurfaces. This gives rise to an {\it{associated curvature function}} $F$, acting on diagonalisable endomorphisms $A$ of an arbitrary real vector space $V$ via
\eq{F(A)=f(\mrm{EV}(A)),}
where $\mrm{EV}(A)$ is the unordered $n$-tuple of eigenvalues of $A$.

 However, when using this definition, $F$ is not defined on the whole space of endomorphisms, but only on the diagonalisable operators. Hence it appears reasonable to view $F$ as defined on bilinear forms,
\eq{\hat{F}(g,h):=F\br{\fr 12 g^{ik}(h_{kj}+h_{jk})}}
for all positive definite $g=(g_{ij})$ and all bilinear forms $h=(h_{ij})\in T^{0,2}_p M$.
Then 
\eq{\hat F^{ij}=\fr{\del F}{\del h_{ij}}}
is a $(2,0)$-tensor and we also write
\eq{\hat F^{ij,kl}=\fr{\del F}{\del h_{ij}\del h_{kl}}.}
Furthermore, if $F=F(\ka_i)$ is strictly monotone, then $\hat{F}^{ij}$ is strictly elliptic. If $F$ is concave, then
\eq{\hat F^{ij,kl}\eta_{ij}\eta_{kl}\leq 0}
for all symmetric $(\eta_{ij})$. We refer to \cite{Andrews:/2007}, \cite[Ch.~2]{Gerhardt:/2006} and \cite{Scheuer:03/2017} for more details on curvature functions. 

Furthermore we will abuse notation and also write $F$ for $\hat{F},$ since no confusion will be possible. E.g., when writing $F^{ij}$, we can only mean $\hat{F}^{ij}$, since there are two contravariant indices. 

We will also use the special curvature functions $H_k$, associated to the $k$-th normalised elementary symmetric polynomial $\s_{k}$ defined on $\Gamma_k$, the connected component of $\{\s_k>0\}$ which contains the point $(1,\dots,1)$.

\subsection{Evolution equations}

 The following evolution equations for \eqref{ICF} are well known and can be found in several places, for example in \cite[Sec.~2.3, Sec.~2.4]{Gerhardt:/2006}. Note that, compared to this reference, we use a different convention on the Riemann tensor.

\begin{lemma}\label{Ev}
Denote $\cF=-F^{-p}$. Along \eqref{ICF} there hold:
\begin{enumerate}[(i)]
\item The induced metric $g$ satisfies
\eq{\label{Ev-g}\dot{g}=-2 \cF h.}
\item The normal vector field satisfies
\eq{\label{Ev-nu}\fr{\bar\n}{dt}\nu=\grad \cF,}
where $\fr{\bar\n}{dt}$ is the covariant time derivative along the curve $x(\cdot,\xi)$ for fixed $\xi\in M$.
\item The second fundamental form satisfies
\eq{\label{Ev-h-1} \dot{h}_{ij}=\cF_{;ij}- \cF h_{ik}h^{k}_{j}+ \cF\ov{\Rm}(x_{;i},\nu,\nu,x_{;j}).}
\item The flow speed $\cF$ satisfies
\eq{\label{Ev-F}\dot{\cF}-\cF^{ij}F_{;ij}&= \cF^{ij}h_{ik}h^{k}_{j}\cF+ \cF^{ij}\ov{\Rm}(x_{;i},\nu,\nu,x_{;j})\cF.}

\end{enumerate}
\end{lemma}

\begin{lemma}\label{Ev-W-2}
Under the flow \eqref{ICF} with $\cF=-F^{-p}$ the second fundamental form evolves by

\eq{\dot{h}_{j}^{i}-\cF^{kl}h^{i}_{j;kl}&=\cF^{kl,rs}{h_{kl;}}^{i}h_{rs;j}+\cF^{kl}h_{rk}h^{r}_{l}h^{i}_{j}-(\cF^{kl}h_{kl}-\cF)h^{i}_{r}h_{j}^{r}\\
				&\hp{=}+\cF^{kl}\bar{R}_{\al\be\g\de}\br{x^{\al}_{;l}x^{\be}_{;j}x^{\g}_{;k}x^{\de}_{;m}h^{im}+x^{\al}_{;l}x^{\be}_{;r}x^{\g}_{;k}x^{\de}_{;m}h^m_jg^{ri}}\\
				&\hp{=}+2\cF^{kl}\bar{R}_{\al\be\g\de}x^{\al}_{;r}x^{\be}_{;m}x^{\g}_{;k}x^{\de}_{;j}h^m_lg^{ri}-\cF^{kl}h_{kl}\bar{R}_{\al\be\g\de}x^{\al}_{;r}\nu^{\be}\nu^{\g}x^{\de}_{;j}g^{ri}\\
				&\hp{=}+ \cF\bar{R}_{\al\be\g\de}x^{\al}_{;r}\nu^{\be}\nu^{\g}x^{\de}_{j}g^{ri}+\cF^{kl}\bar{R}_{\al\be\g\de}x^{\al}_{;k}\nu^{\be}\nu^{\g}x^{\de}_{;l}h^{i}_{j}\\
				&\hp{=}+\cF^{kl}\bar{R}_{\al\be\g\de;\e}\nu^{\al}x^{\be}_{;k}x^{\g}_{;r}x^{\de}_{;l}x^{\e}_{;j}g^{ri}+\cF^{kl}\bar{R}_{\al\be\g\de;\e}\nu^{\al}x^{\be}_{;r}x^{\g}_{;j}x^{\de}_{;k}x^{\e}_{;l}g^{ri}.}
\end{lemma}

\pf{
Basically this is \cite[Lemma~2.4.1]{Gerhardt:/2006}. For convenience we deduce it again, since the proof in that reference is a little rough and we use another convention for the Riemann tensor.
 There hold
\eq{\cF_{;i}=\cF^{kl}h_{kl;i}}
and 
\eq{\label{Ev-h-2-a}\cF_{;ij}&=\cF^{kl,rs}h_{kl;i}h_{rs;j}+\cF^{kl}h_{kl;ij}.}
				
We differentiate the Codazzi equation \eqref{Codazzi-1} to replace the second term on the right hand side. First we differentiate the Codazzi equation with respect to $\del_{j}$, then use the Ricci identities and then differentiate the Codazzi equation with respect to $\del_{l}$. We also use the Weingarten equation \eqref{Weingarten} and the Gauss equation \eqref{GE}.
\eq{
h_{kl;ij}&=h_{ki;lj}-\br{\bar{R}_{\al\be\g\de}\nu^{\al}x^{\be}_{;k}x^{\g}_{;l}x^{\de}_{;i}}_{;j}\\
                &=h_{ki;jl}+{R_{ljk}}^ah_{ai}+{R_{lji}}^a h_{ka}-\br{\bar{R}_{\al\be\g\de}\nu^{\al}x^{\be}_{;k}x^{\g}_{;l}x^{\de}_{;i}}_{;j}\\
                &={R_{ljk}}^ah_{ai}+{R_{lji}}^a h_{ka}-\br{\bar{R}_{\al\be\g\de}\nu^{\al}x^{\be}_{;k}x^{\g}_{;l}x^{\de}_{;i}}_{;j}\\
                &\hp{=}+h_{ij;kl}-\br{\bar{R}_{\al\be\g\de}\nu^{\al}x^{\be}_{;i}x^{\g}_{;k}x^{\de}_{;j}}_{;l}\\
                         &=h_{ij;kl}+(h_{la}h_{jk}-h_{lk}h_{ja}+\bar R_{\al\be\g\de}x^{\al}_{;l}x^{\be}_{;j}x^{\g}_{;k}x^{\de}_{;a})h_{i}^a\\  
                     &\hp{=}+(h_{la}h_{ji}-h_{li}h_{ja}+\bar R_{\al\be\g\de}x^{\al}_{;l}x^{\be}_{;j}x^{\g}_{;i}x^{\de}_{;a})h_{k}^a\\
             &\hp{=}-\bar{R}_{\al\be\g\de;\e}\nu^{\al}x^{\be}_{;k}x^{\g}_{;l}x^{\de}_{;i}x^{\e}_{;j}-\bar{R}_{\al\be\g\de}x^{\al}_{;m}x^{\be}_{;k}x^{\g}_{;l}x^{\de}_{;i}h^m_j+\bar{R}_{\al\be\g\de}\nu^{\al}x^{\be}_{;k}\nu^{\g}x^{\de}_{;i}h_{lj}\\
             &\hp{=}+\bar{R}_{\al\be\g\de}\nu^{\al}x^{\be}_{;k}x^{\g}_{;l}\nu^{\de}h_{ij} -\bar{R}_{\al\be\g\de;\e}\nu^{\al}x^{\be}_{;i}x^{\g}_{;k}x^{\de}_{;j}x^{\e}_{;l}-\bar{R}_{\al\be\g\de}x^{\al}_{;m}x^{\be}_{;i}x^{\g}_{;k}x^{\de}_{;j}h^m_l \\
             &\hp{=}+h_{kl}\bar{R}_{\al\be\g\de}\nu^{\al}x^{\be}_{;i}\nu^{\g}x^{\de}_{;j}+\bar{R}_{\al\be\g\de}\nu^{\al}x^{\be}_{;i}x^{\g}_{;k}\nu^{\de}h_{jl}.
             }

Recall that $h$ satisfies \eqref{Ev-h-1}:
\eq{ \dot{h}_{ij}=\cF_{;ij}- \cF h_{ik}h^{k}_{j}+ \cF\bar{R}_{\al\be\g\de}x^{\al}_{i}\nu^{\be}\nu^{\g}x^{\de}_{j}}
 and hence
 \eq{\dot{h}_{ij}-\cF^{kl}h_{ij;kl}&=\cF^{kl,rs}h_{kl;i}h_{rs;j}-\cF^{kl}h_{kl}h^{r}_{i}h_{rj}+\cF^{kl}h_{rk}h^{r}_{l}h_{ij}\\
				&\hp{=}+\cF^{kl}\bar{R}_{\al\be\g\de}\br{x^{\al}_{;l}x^{\be}_{;j}x^{\g}_{;k}x^{\de}_{;m}h^m_i+x^{\al}_{;l}x^{\be}_{;i}x^{\g}_{;k}x^{\de}_{;m}h^m_j}\\
				&\hp{=}+\cF^{kl}h_{kl}\bar{R}_{\al\be\g\de}\nu^{\al}x^{\be}_{\ ;i}\nu^{\g}x^{\de}_{\ ;j}+\cF^{kl}\bar{R}_{\al\be\g\de}\nu^{\al}x^{\be}_{;k}x^{\g}_{;l}\nu^{\de}h_{ij}\\
				&\hp{=}+2\cF^{kl}\bar{R}_{\al\be\g\de}x^{\al}_{;i}x^{\be}_{;m}x^{\g}_{;k}x^{\de}_{;j}h^m_l-\cF^{kl}\bar{R}_{\al\be\g\de;\e}\nu^{\al}x^{\be}_{;k}x^{\g}_{;l}x^{\de}_{;i}x^{\e}_{;j}\\
				&\hp{=}-\cF^{kl}\bar{R}_{\al\be\g\de;\e}\nu^{\al}x^{\be}_{;i}x^{\g}_{;k}x^{\de}_{;j}x^{\e}_{;l}- \cF h_{ik}h^{k}_{j}+ \cF\bar{R}_{\al\be\g\de}x^{\al}_{i}\nu^{\be}\nu^{\g}x^{\de}_{j}.}

The result follows after reverting to the mixed representation.
}

\subsubsection*{Graphical hypersurfaces}

Given the flow \eqref{ICF} of graphs 
\eq{M_{t}=\{(u(t,y(t,\xi)),y(t,\xi))\cn \xi\in M\}}
 in a warped product with metric of the form \eqref{Warped}, we first of all deduce from \eqref{graph-h} that
 \eq{\label{Ev-u}\dot{u}-\cF^{ij}u_{;ij}&=\fr{p+1}{F^{p}}v^{-1}-\fr{\vt'}{\vt}\fr{p}{F^{p+1}}F^{ij}g_{ij}+\fr{\vt'}{\vt}\fr{p}{F^{p+1}}F^{ij}u_{;i}u_{;j}.}

 Now we deduce the evolution of the quantity
 \eq{w=\fr{1}{\vt^{2}(u)}|du|^{2}_{\s}=|d\p|_{\s}^{2},}
 where $\p$ was defined in \eqref{varphi}.
 The function $\p$ is better suited to these estimates than $u$ itself, since the representation of the second fundamental form is simpler and so the differentiation of the speed $\cF$ is easier to perform. This trick was also used in \cite{Gerhardt:/1990}, \cite{Urbas:/1990} and in subsequent treatments of graphical expanding flows.  Note that $\p$ satisfies
\eq{\label{dotphi}\del_{t}{\p}=-\cF s^{-1},}
where $s$ is the support function defined in \eqref{support}.
In the next lemma we derive the evolution equation for $w$. We simplify notation: Putting lower indices to a function means covariant differentiation with respect to $\s$. 

\begin{lemma}\label{Ev-Gradphi}
Under the flow \eqref{dotphi} in a warped product of the form \eqref{Warped} the gradient function 
\eq{|\hat\n\p|^{2}_{\s}=\p_{i}\p^{i}}
satisfies
\eq{&\br{\fr{d}{dt}-\fr{1}{\vt^{2}}\cF^{k}_{l}\ti{g}^{lr}\hat\n_{kr}}|\hat\n\p|^{2}\\			=~&2\cF \fr{s_{i}}{s^{2}}\p^{i}-2\cF^{k}_{l}h^{l}_{k}\fr{s_{i}}{s^{2}}\p^{i}+4\vt'\cF^{k}_{l}h^{l}_{k}s^{-1}|\hat\n\p|^{2}-2\fr{\vt''}{\vt}\cF^{k}_{k}|\hat\n\p|^{2}\\
		-&\fr{1}{2v^{2}\vt^{2}}\cF^{k}_{l}\s^{lm}|\hat\n\p|^{2}_{m}|\hat\n\p|^{2}_{k}-\fr{1}{v^{2}\vt^{2}}\cF^{k}_{l}\p^{l}|\hat\n\p|^{2}_{r}\p^{r}_{k}\\
		+&\fr{1}{v^{4}\vt^{2}}\cF^{k}_{l}\p^{l}|\hat\n\p|^{2}_{k}|\hat\n\p|^{2}_{i}\p^{i}-\fr{2}{\vt^{2}}\cF^{k}_{l}\ti{g}^{lr}\p_{ir}\p^{i}_{k}-\fr{2}{\vt^{2}}\cF^{k}_{l}\ti{g}^{lr}\hat R_{ikrm}\p^{i}\p^{m}.}
\end{lemma}

\pf{From \eqref{dotphi} we get
\eq{\fr{d}{dt}|\hat\n\p|^{2}&=2\dot{\p}_{i}\p^{i}=\fr{2}{s^{2}}\cF s_{i}\p^{i}-2\cF^{k}_{l}\hat{\n}_{i}h^{l}_{k}\p^{i}s^{-1}.}

Due to \eqref{Wphi} there holds
\eq{\hat{\n}_{i}h^{l}_{k}&=-\fr{v_{i}\vt+v\vt'\vt\p_{i}}{v^{2}\vt^{2}}\br{\vt'\de^{l}_{k}-\ti g^{lr}\p_{rk}}+\fr{1}{v\vt}\br{\vt''\vt\p_{i}\de^{l}_{k}-\hat\n_{i}\ti g^{lr}\p_{rk}-\ti g^{lr}\p_{rki}}\\
				&=-\fr{v_{i}}{v}h^{l}_{k}-\vt'\p_{i}h^{l}_{k}+\fr{\vt''}{v}\de^{l}_{k}\p_{i}+\fr{\p^{l}_{i}\p^{r}+\p^{l}\p^{r}_{i}}{v^{3}\vt}\p_{rk}-\fr{2}{v^{4}\vt}v_{i}\p^{l}\p^{r}\p_{rk}\\
				&\hp{=}-\fr{1}{v\vt}\ti g^{lr}\p_{rki}\\
				&=\fr{s_{i}}{s}h^{l}_{k}-2\vt'\p_{i}h^{l}_{k}+\fr{\vt''}{v}\de^{l}_{k}\p_{i}+\fr{\p^{l}_{i}\p^{r}+\p^{l}\p^{r}_{i}}{v^{3}\vt}\p_{rk}\\
				&\hp{=}-\fr{1}{v^{5}\vt}|\hat\n\p|^{2}_{i}\p^{l}\p^{r}\p_{rk}-\fr{1}{v\vt}\ti g^{lr}\p_{rik}+\fr{1}{v\vt}\ti{g}^{lr}{\hat R_{ikr}}^{m}\p_{m},}
where we used the definition of the Riemann tensor of $\s$.
Using
\eq{|\hat\n\p|^{2}_{rk}=2\p_{irk}\p^{i}+2\p_{ir}\p^{i}_{k},}
we combine these two equalities to get
\eq{&\br{\fr{d}{dt}-\fr{1}{\vt^{2}}\cF^{k}_{l}\ti{g}^{lr}\hat\n_{kr}}|\hat\n\p|^{2}\\			=~&2\cF \fr{s_{i}}{s^{2}}\p^{i}-2\cF^{k}_{l}h^{l}_{k}\fr{s_{i}}{s^{2}}\p^{i}+4\vt'\cF^{k}_{l}h^{l}_{k}s^{-1}|\hat\n\p|^{2}-2\fr{\vt''}{\vt}F^{k}_{k}|\hat\n\p|^{2}\\
		-&\fr{2}{v^{2}\vt^{2}}\cF^{k}_{l}(\p^{l}_{i}\p^{r}+\p^{l}\p^{r}_{i})\p^{i}\p_{rk}+\fr{2}{v^{4}\vt^{2}}\cF^{k}_{l}\p^{l}\p^{r}\p_{rk}|\hat\n\p|^{2}_{i}\p^{i}\\
		-&\fr{2}{\vt^{2}}\cF^{k}_{l}\ti{g}^{lr}\p_{ir}\p^{i}_{k}-\fr{2}{\vt^{2}}\cF^{k}_{l}\ti{g}^{lr}\hat R_{ikr}^{m}\p^{i}\p_{m}}
and hence the result.
}

The support function satisfies the following evolution.
\begin{lemma}
Along \eqref{ICF} in a warped product with metric \eqref{Warped}, the support function
\eq{s=\vt(u)\bar{g}(\del_{r},\nu)}
satisfies 
\eq{\label{Ev-s}\dot{s}-\cF^{ij}s_{;ij}&=\cF^{ij}h_{ik}h^{k}_{j}s-\vt'\fr{p-1}{F^{p}}+\bar g(\vt\del_{r},\n\cF)-\cF^{ij}(\bar g(\vt\del_{r},x_{;k}h^{k}_{i;j})).}
\end{lemma}

\pf{
The vector field $\vt\del_{r}$ is conformal,
\eq{\bar\n_{\bar X}(\vt\del_{r})=\vt'\bar X\q\fa \bar X\in T^{1,0}(N).}
Hence
\eq{\label{Ev-s-1}\dot{s}=\bar{g}(\vt'\dot{x},\nu)+\bar{g}(\vt\del_{r},\bar\n_{\dot{x}}\nu)=-\vt'\cF+\bar{g}(\vt\del{r},\n\cF),}
\eq{\label{grad-s}Xs=\bar{g}(\vt\del_{r},\cW(X))}
and 
\eq{\label{Hess-s}\n^{2}s(X,Y)&=Y(Xs)-(\n_{Y}X)s\\
			&=\vt'h(X,Y)-h(X,\cW(Y))s+\bar{g}(\vt\del_{r},\n_{Y}\cW(X)).}
The result follows from combining these equalities.
}

We will also make use of the evolution of $\dot{\p}$. This method was used in \cite[Prop.~3.4]{BrendleHungWang:01/2016} and \cite[Lemma~3.5]{Scheuer:01/2017}.

\begin{lemma}
Under the flow \eqref{dotphi} in a warped product of the form \eqref{Warped} the speed $\dot\p$ satisfies
\eq{\label{L-dotphi}\del_{t}\dot\p-\fr{\del\dot{\p}}{\del\p_{ij}}\dot{\p}_{ij}-\fr{\del \dot{\p}}{\del\p_{i}}\dot{\p}_{i}=\fr{\vt'}{\vt}v\cF^{i}_{j}h^{j}_{i}\dot{\p}-\fr{\vt''}{\vt}\cF^{k}_{k}\dot{\p}+\fr{\vt'}{\vt}v\cF \dot\p.}
\end{lemma}

\pf{Differentiating 
\eq{\dot{\p}=-\cF s^{-1}}
gives
\eq{\del_{t}\dot\p-\fr{\del \dot{\p}}{\del \p_{ij}}\dot{\p}_{ij}-\fr{\del \dot{\p}}{\del \p_{i}}\dot{\p}_{i}&=\fr{\del \dot{\p}}{\del \p}\dot{\p}\\
			&=-\cF^{j}_{i}\fr{h^{i}_{j}}{\del\p} s^{-1}\dot\p+s^{-1}\cF\vt' \dot\p.}
From \eqref{Wphi} we get
\eq{\fr{\del h^{i}_{j}}{\del\p}&=-\fr{\vt'}{v\vt}\br{\vt'\de^{i}_{j}-\ti{g}^{ik}\p_{kj}}+\fr{\vt''}{v}\de^{i}_{j}}
and inserting this gives the result.
}

\section{Long-time existence}\label{LTE}
\subsection{Barriers}

\begin{lemma}\label{ICF-Barriers}
Let $\vt\in C^{2}((R_{0},\8))$ with $\vt'>0$ and $\vt''\geq 0$, $r_{0}>R_{0}$ and $0<p\leq 1$.
Let $r(t,r_{0})$ be the unique solution of the initial value problem
\begin{align}\dot{r}&=\fr{\vt^{p}(r)}{n^{p}\vt'^{p}(r)}\label{ODE}\\
		r(0)&=r_{0}.\end{align}
Then $r$ is defined for all times and
\eq{r(t,r_{0})\ra \8,\q t\ra \8.} 
Consequently, for $x_{0}$ as in \cref{ICF-Main} with associated graph function $u_{0}$, we have
\eq{\inf_{M}u(t,\cdot)\ra \8,\q t\ra \8,}
provided the flow \eqref{ICF-Ini} exists for all times.
\end{lemma}

\pf{
Due to \eqref{ODE} we have
\eq{\dot{r}\leq 1+\fr{\vt(r)}{\vt'(r)},}
where the right hand side grows at most linearly in $r$ due to $\vt''\geq 0$. Hence $r$ is defined for all times.
Suppose $r$ does not converge to infinity. Due to its monotonicity it converges to some $r_{1}<\8$. From
\eq{\fr{\vt}{\vt'}(r)>0\q\fa r\in[r_{0},r_{1}]}
we obtain $\dot{r}\geq c>0$ and reach a contradiction. The second claim follows from the maximum principle which gives
\eq{r(t,\inf u_{0})\leq u(t,\cdot)\leq r(t,\sup u_{0}).}
}

\subsection{Gradient estimates}

Let us first prove some rough gradient estimates which will suffice to get the long-time existence.
In the a priori estimates that appear in the rest of the paper, generic constants will be allowed to depend on the data of the problem, namely $N,p,M_{0}$ unless otherwise stated.

 First we need a bound on $F$ from below:

\begin{lemma}\label{ICF-1/F-bound}
Under the assumptions of \cref{ICF-Main}~(i), along \eqref{ICF} the spatial maxima of the quantity
\eq{\dot{\p}=\fr{1}{F^{p}}\fr{v}{\vt}}
are non-increasing.
\end{lemma}

\pf{
According to \eqref{L-dotphi}
we have
\eq{\del_{t}\dot{\p}-\fr{\del\dot{\p}}{\del \p_{ij}}\dot{\p}_{ij}-\fr{\del \dot{\p}}{\del\p_{i}}\dot{\p}_{i}=\fr{(p-1)\vt'}{\vt F^{p}}v\dot{\p}-\fr{\vt''p}{\vt F^{p+1}}F^{i}_{i}\dot{\p}\leq 0.}

The result follows from the maximum principle.
}

Now we prove some very general gradient estimates for inverse curvature flows in warped products. We use the notation from \cref{Ev-Gradphi}.

\begin{lemma}\label{ICF-grad-bound}
Under the assumptions of \cref{ICF-Main}~(i), along \eqref{ICF} the function
$|\hat\n\p|^{2}$
is bounded on every finite time interval. Furthermore, under the assumptions of \cref{ICF-Main}~(ii), there exists a positive constant $\g$, such that the spatial maxima of
\eq{\hat{z}=|\hat\n\p|^{2}\vt^{\g}}
are non-increasing, provided $p<1$, and such that the spatial maxima of 
\eq{\tilde{z}=|\hat\n\p|^{2}\vt'^{\g}}
are non-increasing, regardless the value of $0<p\leq 1$.
\end{lemma}

\pf{
We want to calculate the evolution equations of $\hat z$ and $\ti{z}$. Hence we need one for $u$, which makes use of the parabolic operator with respect to the metric $\s$. Note that in \eqref{Ev-u} we use covariant derivatives of the metric induced by $u$, hence we need to rewrite this. The covariant derivatives with respect to $\s$ and $g$ are related by
\eq{\hat\n^{2}u&=v^{2}\n^{2}u+\fr{\vt'}{\vt}\br{2du\otimes du-\vt^{2}(v^{2}-1)\s}\\
			&= -vh+v^{2}\fr{\vt'}{\vt}g-v^{2}\fr{\vt'}{\vt}du\otimes du+2\fr{\vt'}{\vt}du\otimes du-\vt'\vt(v^{2}-1)\s, }
cf. \cite[equ.~(71)]{Scheuer:01/2017} and \eqref{graph-h}.
We obtain
\eq{\br{\del_{t}-\fr{1}{\vt^{2}}\fr{p}{F^{p+1}}F^{k}_{l}\ti{g}^{lr}\hat\n^{2}_{kr}}u&=\fr{v}{F^{p}}+\fr{pv}{F^{p}}-\fr{p}{F^{p+1}}\fr{\vt'}{\vt}F^{k}_{k}v^{2}\\
	&\hp{=}+\fr{p}{F^{p+1}}\fr{\vt'}{\vt}F^{kr}u_{;k}u_{;r}v^{2}-\fr{2p}{F^{p+1}}\fr{\vt'}{\vt}F^{kr}u_{;k}u_{;r}\\
	&\hp{=}+\fr{p}{F^{p+1}}\fr{\vt'}{\vt}F^{kr}(g_{kr}-u_{;k}u_{;r})(v^{2}-1)\\
	&=\fr{p+1}{F^{p}}v-\fr{p}{F^{p+1}}\fr{\vt'}{\vt}F^{k}_{k}-\fr{p}{F^{p+1}}\fr{\vt'}{\vt}F^{kr}u_{;k}u_{;r}.}
	
Now we use
\eq{\fr{s_{i}}{s}=\vt'\p_{i}-\fr{v_{i}}{v}=\vt'\p_{i}-\fr{1}{2}\fr{|\hat\n\p|^{2}_{i}}{v^{2}}}
to deduce from \cref{Ev-Gradphi}:
\eq{&\br{\del_{t}-\fr{1}{\vt^{2}}\fr{p}{F^{p+1}}F^{k}_{l}\ti{g}^{lr}\hat\n_{kr}}|\hat\n\p|^{2}\\
=~&-\fr{2}{F^{p}}\fr{\vt'}{\vt}v|\hat\n\p|^{2}+\fr{1}{\vt vF^{p}}|\hat\n\p|^{2}_{i}\p^{i}-\fr{2p}{F^{p}}\fr{\vt'}{\vt}v|\hat\n\p|^{2}+\fr{p}{\vt vF^{p}}|\hat\n\p|^{2}_{i}\p^{i}\\
\hp{=}+&\fr{4p}{F^{p}}\fr{\vt'}{\vt}v|\hat\n\p|^{2}-\fr{2p}{F^{p+1}}\fr{\vt''}{\vt}F^{k}_{k}|\hat\n\p|^{2}-\fr{1}{2v^{2}\vt^{2}}\fr{p}{F^{p+1}}F^{k}_{l}\s^{lm}|\hat\n\p|^{2}_{m}|\hat\n\p|^{2}_{k}\\
\hp{=}-&\fr{1}{v^{2}\vt^{2}}\fr{p}{F^{p+1}}F^{k}_{l}\p^{l}|\hat\n\p|^{2}_{r}\p^{r}_{k}+\fr{1}{v^{4}\vt^{2}}\fr{p}{F^{p+1}}F^{k}_{l}\p^{l}|\hat\n\p|^{2}_{k}|\hat\n\p|^{2}_{i}\p^{i}\\
\hp{=}-&\fr{2}{\vt^{2}}\fr{p}{F^{p+1}}F^{k}_{l}\ti{g}^{lr}\p_{ir}\p^{i}_{k}-\fr{2}{\vt^{2}}\fr{p}{F^{p+1}}F^{k}_{l}\ti{g}^{lr}{\hat{R}_{ikr}}^{m}\p^{i}\p_{m}\\
=~&\fr{2(p-1)}{F^{p}}\fr{\vt'}{\vt}v|\hat\n\p|^{2}-\fr{2p}{F^{p+1}}\fr{\vt''}{\vt}F^{k}_{k}|\hat\n\p|^{2}-\fr{2}{\vt^{2}}\fr{p}{F^{p+1}}F^{k}_{l}\ti{g}^{lr}\hat{R}_{ikr}^{m}\p^{i}\p_{m}\\
\hp{=}+&\fr{p+1}{F^{p}}\fr{1}{v\vt}|\hat\n\p|^{2}_{i}\p^{i}-\fr{2p}{F^{p+1}}\fr{1}{\vt^{2}}F^{k}_{l}\ti{g}^{lr}\p_{ir}\p^{i}_{k}-\fr{1}{2v^{2}\vt^{2}}\fr{p}{F^{p+1}}F^{k}_{l}\s^{lm}|\hat\n\p|^{2}_{k}|\hat\n\p|^{2}_{m}\\
\hp{=}+&\fr{1}{v^{4}\vt^{2}}\fr{p}{F^{p+1}}F^{k}_{l}\p^{l}|\hat\n\p|^{2}_{k}|\hat\n\p|^{2}_{i}\p^{i}-\fr{1}{v^{2}\vt^{2}}\fr{p}{F^{p+1}}F^{k}_{l}\p^{l}|\hat\n\p|^{2}_{r}\p^{r}_{k}.
}

Now first generally put
\eq{z=f(u)|\hat\n\p|^{2}.}
With the help of the previous calculations we get at a maximal point of $z$, where
\eq{|\hat\n\p|^{2}_{i}=-\fr{f'}{f}\vt|\hat\n\p|^{2}\p_{i},}
\eq{\label{ICF-grad-bound-1}&\cL z\equiv\br{\del_{t}-\fr{1}{\vt^{2}}\fr{p}{F^{p+1}}F^{k}_{l}\ti{g}^{lr}\hat\n_{kr}}z\\
=~&f\cL|\hat\n\p|^{2}+z\fr{f'}{f}\cL u-\fr{f''}{f}\fr{p}{\vt^{2}F^{p+1}}F^{k}_{l}\ti{g}^{lr}u_{k}u_{r}z-\fr{2p}{\vt^{2}F^{p+1}}F^{k}_{l}\ti{g}^{lr}|\hat\n\p|^{2}_{k}f_{r}\\
=~&\fr{2(p-1)}{F^{p}}\fr{\vt'}{\vt}vz-\fr{2p}{F^{p+1}}\fr{\vt''}{\vt}F^{k}_{k}z-\fr{2f}{\vt^{2}}\fr{p}{F^{p+1}}F^{k}_{l}\ti{g}^{lr}\hat{R}_{ikr}^{m}\p^{i}\p_{m}\\
\hp{=}-&\fr{p+1}{F^{p}}\fr{f'}{f}\fr{1}{v}|\hat\n\p|^{2}z-\fr{2p}{F^{p+1}}\fr{1}{\vt^{2}}F^{k}_{l}\ti{g}^{lr}\p_{ir}\p^{i}_{k}f\\
\hp{=}-&\fr{1}{2v^{2}}\fr{f'^{2}}{f^{2}}\fr{p}{F^{p+1}}F^{k}_{l}\p^{l}\p_{k}|\hat\n\p|^{2}z+\fr{1}{v^{4}}\fr{f'^{2}}{f^{2}}\fr{p}{F^{p+1}}F^{k}_{l}\p^{l}\p_{k}|\hat\n\p|^{4}z\\
\hp{=}-&\fr{1}{2v^{2}}\fr{f'^{2}}{f^{2}}\fr{p}{F^{p+1}}F^{k}_{l}\p^{l}\p_{k}|\hat\n\p|^{2}z+\fr{p+1}{F^{p}}\fr{f'}{f}vz-\fr{p}{F^{p+1}}\fr{f'}{f}\fr{\vt'}{\vt}F^{k}_{k}z\\
\hp{=}-&\fr{p}{F^{p+1}}\fr{f'}{f}\fr{\vt'}{\vt}F^{kr}u_{;k}u_{;r}z+\fr{2p}{F^{p+1}}\fr{f'^{2}}{f^{2}}F^{k}_{l}\ti{g}^{lr}\p_{k}\p_{r}z-\fr{p}{F^{p+1}}\fr{f''}{f}F^{kr}u_{k}u_{r}z\\
=~&\fr{2(p-1)}{F^{p}}\fr{\vt'}{\vt}vz-\fr{2p}{F^{p+1}}\fr{\vt''}{\vt}F^{k}_{k}z-\fr{2f}{\vt^{2}}\fr{p}{F^{p+1}}F^{k}_{l}\ti{g}^{lr}\hat{R}_{ikr}^{m}\p^{i}\p_{m}\\
\hp{=}+&\fr{p+1}{F^{p}}\fr{f'}{f}\fr{1}{v}z-\fr{2p}{F^{p+1}}\fr{1}{\vt^{2}}F^{k}_{l}\ti{g}^{lr}\p_{ir}\p^{i}_{k}f-\fr{f'^{2}}{f^{2}}\fr{p}{F^{p+1}}F^{k}_{l}\p^{l}\p_{k}\fr{|\hat\n\p|^{2}}{v^{4}}z\\
\hp{=}-&\fr{p}{F^{p+1}}\fr{f'}{f}\fr{\vt'}{\vt}F^{k}_{k}z+\fr{p}{F^{p+1}}F^{kr}u_{k}u_{r}z\br{2\fr{f'^{2}}{f^{2}}-\fr{f''}{f}-\fr{f'}{f}\fr{\vt'}{\vt}}.}
Employing \cref{ICF-1/F-bound} with $f=1$ we see $\cL z\leq cz,$
where on each finite interval this constant is bounded. We also used
\eq{F=n\fr{H_{k+1}}{H_{k}}\q\Ra \q F^{k}_{k}\leq c,}
cf. \cite[Lemma~2.7]{Lu:09/2016}. Hence we obtain the first claim.
Under the assumptions of \cref{ICF-Main} (ii), we obtain that $\max \hat z$ is decreasing if $p<1$ and $f(u)=\vt^{\g}(u)$
for a small $\g>0$. If $p=1$, the same is true with $\g=0$.
Due to \eqref{Wphi} we have
\eq{F=F^{k}_{l}h^{l}_{k}=\fr{\vt'}{v\vt}F^{k}_{k}-\fr{1}{v\vt}F^{k}_{l}\ti{g}^{lr}\p_{rk}}
 and hence
 \eq{&\fr{p+1}{F^{p}}\fr{f'}{f}\fr{1}{v}z-\fr{2p}{F^{p+1}}\fr{1}{\vt^{2}}F^{k}_{l}\ti{g}^{lr}\p_{ir}\p^{i}_{k}f\\
 	=~&\fr{p+1}{F^{p+1}}\fr{f'}{f}\fr{\vt'}{\vt}\fr{1}{v^{2}}F^{k}_{k}z-\fr{p}{F^{p+1}}\fr{p+1}{p}\fr{f'}{f}\fr{1}{v^{2}\vt}F^{k}_{l}\ti{g}^{lr}\p_{rk}z-\fr{2p}{F^{p+1}}\fr{1}{\vt^{2}}F^{k}_{l}\ti{g}^{lr}\p_{ir}\p^{i}_{k}f\\
	=~&\fr{p+1}{F^{p+1}}\fr{f'}{f}\fr{\vt'}{\vt}\fr{1}{v^{2}}F^{k}_{k}z+\fr{p}{F^{p+1}}\fr{(p+1)^{2}}{8p}\fr{f'^{2}}{f^{2}}\fr{|\hat\n\p|^{2}}{v^{4}}F^{k}_{l}\ti{g}^{lr}\s_{rk}z\\
	\hp{=}-&\fr{p}{F^{p+1}}\fr{f}{2}F^{k}_{l}\ti{g}^{lr}\br{\fr{p+1}{2p}\fr{f'}{f^{2}}\fr{1}{v^{2}}z\s_{ir}+\fr{2}{\vt}\p_{ir}}\br{\fr{p+1}{2p}\fr{f'}{f^{2}}\fr{1}{v^{2}}z\de^{i}_{k}+\fr{2}{\vt}\p^{i}_{k}}.}
Hence at a maximal point of $z$ we get for $f=\vt'^{\g}$
\eq{\label{ICF-grad-bound-2}\cL z&\leq -\fr{2p}{F^{p+1}}\br{\fr{\vt''}{\vt}+\fr{1}{2}\fr{f'}{f}\fr{\vt'}{\vt}-\fr{(p+1)^{2}}{16p}\fr{f'^{2}}{f^{2}}|\hat\n\p|^{2}-\fr{p+1}{2p}\fr{f'}{f}\fr{\vt'}{\vt}}F^{k}_{k}z\\
&\hp{=}+\fr{2(p-1)}{F^{p}}\fr{\vt'}{\vt}vz-\fr{2f}{\vt^{2}}\fr{p}{F^{p+1}}F^{k}_{l}\ti{g}^{lr}\hat{R}_{ikr}^{m}\p^{i}\p_{m}\\
&\hp{=}+\fr{p}{F^{p+1}}F^{kr}u_{k}u_{r}z\br{2\fr{f'^{2}}{f^{2}}-\fr{f''}{f}-\fr{f'}{f}\fr{\vt'}{\vt}},}
which is negative if $0<\g$ is small enough. Here we also used \cref{ICF-N}.
}

\subsection{Curvature estimates}
We prove that along \eqref{ICF} all principal curvatures are bounded as long as the flow remains in a compact subset of $N$. Due to all previous a priori estimates this will imply uniform $C^{2}$ estimates on each finite time interval, as well as a uniformly elliptic operator $F^{-(p+1)}F^{ij}$. Hence the regularity estimates by Krylov and Safonov apply to get $C^{2,\al}$ estimates. With the linear Schauder estimates we obtain uniform $C^{\8}$-bounds on each finite interval. We may extend the solution beyond any finite $T$, completing the proof of item (i) of \cref{ICF-Main}.

\begin{prop}\label{ICF-k-bound}
Under the assumptions of \cref{ICF-Main}~(i), on every finite interval $[0,T]$ there exists a compact set $\La\sub\G$ such that along the flow \eqref{ICF} the principal curvatures $\ka_{i}$ satisfy
\eq{(\ka_{i})\in \La\q\fa 0\leq t\leq T.}
\end{prop}

\pf{In this proof the generic constant $c$ is allowed to depend on $T$. We proceed similarly as in \cite{Gerhardt:11/2011,ScheuerXia:08/2017}. First we simplify the evolution of the second fundamental form, cf. \cref{Ev-W-2}. We have the following estimate in normal coordinates:
\eq{\dot{h}^{n}_{n}-\cF^{kl}h^{n}_{n;kl}&\leq \cF^{kl,rs}h_{kl;n}{h_{rs;}}^{n}+\fr{p}{F^{p+1}}F^{kl}h_{rk}h^{r}_{l}h^{n}_{n}-\fr{p+1}{F^{p}}(h^{n}_{n})^{2}\\
			&\hp{=}+\fr{c}{F^{p+1}}F^{ij}g_{ij}(h^{n}_{n}+1)+\fr{c}{F^{p}}.}According to \eqref{Ev-s} the support function,
which is bounded from below due to the gradient estimates \cref{ICF-grad-bound},
satisfies
\eq{\label{ICF-Ev-s}\dot{s}-\cF^{kl}s_{;kl}&=\fr{p}{F^{p+1}}F^{kl}h_{rk}h^{r}_{l}s-\fr{(p-1)\vt'}{F^{p}}\\
				&\hp{=}+\fr{p\vt}{F^{p+1}}F^{kl}{u_{;}}^{m}\ov{\Rm}(\nu,x_{;k},x_{;m},x_{;l}).}
Due to a well known trick, e.g. see the proof of \cite[Lemma~4.4]{Gerhardt:11/2011}, it suffices to bound the evolution equation of the function
\eq{w=\log h^{n}_{n}+f(s)+\al u,}
$\al$ to be determined,
at a maximal point of $w$ in which normal coordinates are given,
\eq{g_{ij}=\de_{ij},\q h_{ij}=\ka_{i}\de_{ij},\q \ka_{1}\leq\dots\leq\ka_{n}.}
For small $\be>0$ set
\eq{f(s)=-\log(s-\be).}
Also using \eqref{Ev-u}, we see that $w$ satisfies
\eq{\label{ICF-k-bound-1}\dot{w}-&\cF^{kl}w_{;kl}\leq \fr{p}{F^{p+1}}F^{kl}h_{rk}h^{r}_{l}(1+f's)-\fr{p+1}{F^{p}}h^{n}_{n}\\
				&\hp{=}+\fr{c}{F^{p+1}}F^{k}_{k}(1+(h^{n}_{n})^{-1})+\fr{c}{F^{p}}(1+(h^{n}_{n})^{-1}+\al)\\
				&\hp{=}+\fr{p}{F^{p+1}}F^{kl,rs}h_{kl;n}{h_{rs;}}^{n}(h^{n}_{n})^{-1}+\fr{p}{F^{p+1}}F^{ij}(\log h^{n}_{n})_{;i}(\log h^{n}_{n})_{;j}\\
				&\hp{=}-f''\fr{p}{F^{p+1}}F^{ij}s_{;i}s_{;j}-\al\fr{\vt'}{\vt}\fr{p}{F^{p+1}}F^{ij}(g_{ij}-u_{;i}u_{;j}).}

We employ a trick already used in \cite{Enz:10/2008}.
Due to the concavity of $F$ there holds
\eq{F^{nn}\leq \dots \leq F^{11}\q \text{and}\q F^{kl,rs}\eta_{kl}\eta_{rs}\leq \fr{2}{\ka_{n}-\ka_{1}} \sum_{k=1}^{n}(F^{nn}-F^{kk})\eta_{nk}^{2}}
for all symmetric $(\eta_{ij})$. It is possible to exploit this term in order to estimate \eqref{ICF-k-bound-1}.

{\bf{Case 1:}} $\ka_{1}<-\e_{1}\ka_{n}$, $0<\e_{1}<\fr 12$. There hold
\eq{F^{ij}h_{ik}h^{k}_{j}\geq F^{11}\ka_{1}^{2}\geq \fr{1}{n}F^{ij}g_{ij}\ka_{1}^{2}\geq \fr 1n F^{ij}g_{ij}\e_{1}^{2}\ka_{n}^{2},}
\eq{\fr{p}{F^{p+1}}F^{ij}(\log h^{n}_{n})_{;i}(\log h^{n}_{n})_{;j}
	&=f'^{2}\fr{p}{F^{p+1}}F^{ij}s_{;i}s_{;j}+ f'\fr{2\al p}{F^{p+1}}F^{ij}s_{;i}u_{;j}	\\
	&\hp{=}+\fr{\al^{2}p}{F^{p+1}}F^{ij}u_{;i}u_{;j},}
due to $\n w=0$. Using \eqref{grad-s}, if $\ka_{n}$ is large, \eqref{ICF-k-bound-1} with any $\al$ becomes
\eq{ 0&\leq \fr{p}{nF^{p+1}}F^{ij}g_{ij}\br{\e_{1}^{2}\ka_{n}^{2}(1+f's)+c+c\al|f'|\ka_{n}+c\al^{2}}\\
		&\hp{=}+\fr{p+1}{F^{p}}\br{c+c\al-\ka_{n}}-\fr{p}{F^{p+1}}F^{ij}s_{;i}s_{;j}(f''-f'^{2})\\
		&<0}
and hence $w$ is bounded in this case, due to $1+f's\leq c<0$ and $f''=f'^{2}$.

{\bf{Case 2:}} $\ka_{1}\geq -\e_{1}\ka_{n}.$ Then
\eq{&\fr{2}{\ka_{n}-\ka_{1}}\sum_{k=1}^{n}(F^{nn}-F^{kk})(h_{nk;n})^{2}(h^{n}_{n})^{-1}\\
	\leq~& \fr{2}{1+\e_{1}}\sum_{k=1}^{n}(F^{nn}-F^{kk})(h_{nk;n})^{2}(h^{n}_{n})^{-2}\\
	\leq~&\fr{2}{1+\e_{1}}\sum_{k=1}^{n}(F^{nn}-F^{kk})(h_{nn;k})^{2}(h^{n}_{n})^{-2}+c(\e_{1})\sum_{k=1}^{n}(F^{kk}-F^{nn})\ka_{n}^{-2}\\
	&\hp{=}+\fr{4}{1+\e_{1}}\sum_{k=1}^{n}(F^{nn}-F^{kk})h_{nn;k}\bar{R}_{\al\be\g\de}\nu^{a}x^{\be}_{;n}x^{\g}_{;n}x^{\de}_{;k}(h^{n}_{n})^{-2}\\
	\leq~& \fr{2}{1+2\e_{1}}\sum_{k=1}^{n}(F^{nn}-F^{kk})(h_{nn;k})^{2}(h^{n}_{n})^{-2}+c(\e_{1})\sum_{k=1}^{n}(F^{kk}-F^{nn})\ka_{n}^{-2}, }
where we used the Codazzi equation \eqref{Codazzi-1} and Cauchy-Bunjakowski-Schwarz.
We deduce further:
\eq{&F^{ij}(\log h^{n}_{n})_{;i}(\log h^{n}_{n})_{;j}+\fr{2}{\ka_{n}-\ka_{1}}\sum_{k=1}^{n}(F^{nn}-F^{kk})(h_{nk;n})^{2}(h^{n}_{n})^{-1}\\
	\leq~&\fr{2}{1+2\e_{1}}\sum_{k=1}^{n}F^{nn}(\log h^{n}_{n})_{;k}^{2}-\fr{1-2\e_{1}}{1+2\e_{1}}\sum_{k=1}^{n}F^{kk}(\log h^{n}_{n})_{;k}^{2}+c(\e_{1})F^{ij}g_{ij}\ka_{n}^{-2}\\
	\leq~&\sum_{k=1}^{n}F^{nn}(\log h^{n}_{n})_{;k}^{2}+c(\e_{1})F^{ij}g_{ij}\ka_{n}^{-2}\\
	=~&c(\e_{1})F^{ij}g_{ij}\ka_{n}^{-2}+f'^{2}F^{nn}\|\n s\|^{2}+2\al f'F^{nn}\ip{\n s}{\n u}+\al^{2}F^{nn}\|\n u\|^{2}.}

Hence we can estimate \eqref{ICF-k-bound-1}, using $F^{ij}\bar g_{ij}\geq c_{0}F^{ij}g_{ij},$

\eq{0&\leq\fr{p}{F^{p+1}}F^{nn}\br{\ka_{n}^{2}(1+f's)+\al c\ka_{n}+c\al^{2}}+\fr{p+1}{F^{p}}\br{c+\al-\ka_{n}}\\
		&\hp{=}+\fr{p}{F^{p+1}}F^{ij}g_{ij}\br{\fr{c(\e_{1})}{\ka_{n}^{2}}+c-c_{0}\al\fr{\vt'}{\vt}}+\fr{p}{F^{p+1}}\br{f'^{2}F^{nn}\|\n s\|^{2}-f''F^{ij}s_{;i}s_{;j}}.} 
Due to the barrier estimates, on every finite interval $[0,T]$ there holds $\vt\leq c(T).$
Picking $\al$ large enough, we see that $\ka_{n}$ is bounded on $[0,T]$ and the proof is complete.
}

\section{Asymptotics}\label{Asym}
\subsection{Global bounds}

In order to study the long-time behaviour of \eqref{ICF}, we need to investigate the evolution of the second fundamental form in greater detail. Therefore we need a more detailed version of its evolution equation.

\begin{lemma}\label{ICF-Ev-W}
Along \eqref{ICF} the Weingarten operator evolves according to
\eq{&\dot{h}^{i}_{j}-\fr{p}{F^{p+1}}F^{kl}h^{i}_{j;kl}\\
=~&-\fr{p(p+1)}{F^{p+2}}F_{;j}{F_{;}}^{i}+\fr{p}{F^{p+1}}F^{kl,rs}h_{kl;j}{h_{rs;}}^{i}+\fr{p}{F^{p+1}}F^{kl}h_{rk}h^{r}_{l}h^{i}_{j}-\fr{p+1}{F^{p}}h^{i}_{r}h^{r}_{j}\\
+&~\fr{\vt'^{2}}{\vt^{2}}\fr{p}{F^{p+1}}F^{k}_{k}h^{i}_{j}-\fr{\vt'^{2}}{\vt^{2}}\fr{p-1}{F^{p}}\de^{i}_{j}-\br{\fr{\vt''}{\vt}-\fr{\vt'^{2}}{\vt^{2}}}v^{-2}\fr{p}{F^{p+1}}F^{k}_{k}h^{i}_{j}\\
+&~\br{\fr{\vt''}{\vt}-\fr{\vt'^{2}}{\vt^{2}}}v^{-2}\fr{p+1}{F^{p}}\de^{i}_{j}+\br{\fr{\vt''}{\vt}-\fr{\vt'^{2}}{\vt^{2}}}\fr{p}{F^{p+1}}F^{kl}\br{u_{;k}u_{;l}h^{i}_{j}-2h^{m}_{l}u_{;m}u_{;k}\de^{i}_{j}}\\
+&~\br{\fr{\vt''}{\vt}-\fr{\vt'^{2}}{\vt^{2}}}\fr{1-p}{F^{p}}u_{;j}{u_{;}}^{i}+\br{\fr{\vt''}{\vt}-\fr{\vt'^{2}}{\vt^{2}}}\fr{p}{F^{p+1}}F^{k}_{k}\br{h^{i}_{m}{u_{;}}^{m}u_{;j}+h^{m}_{j}u_{;m}{u_{;}}^{i}}\\
				+&~\br{\fr{\vt''}{\vt}-\fr{\vt'^{2}}{\vt^{2}}}\fr{p}{F^{p+1}}\br{F^{il}u_{;m}(h^{m}_{l}u_{;j}-h^{m}_{j}u_{;l})+F^{k}_{l}h^{l}_{j}{u_{;}}^{i}u_{;k}-F^{l}_{j}h^{im}u_{;l}u_{;m}}\\
	+&~\fr{p}{F^{p+1}}F^{kl}\br{\widetilde{\Rm}(x_{;l},x_{;j},x_{;k},x_{;m})h^{im}+\widetilde{\Rm}(x_{;l},x_{;r},x_{;k},x_{;m})h^{m}_{j}g^{ri}}\\
	+&~\fr{2p}{F^{p+1}}F^{kl}\widetilde{\Rm}(x_{;r},x_{;m},x_{;k},x_{;j})h^{m}_{l}g^{ri}-\fr{p+1}{F^{p}}\widetilde{\Rm}(x_{;r},\nu,\nu,x_{;j})g^{ri}\\
	+&~\fr{p}{F^{p+1}}F^{kl}\widetilde{\Rm}(x_{;k},\nu,\nu,x_{;l})h^{i}_{j}\\
	+&~\fr{p}{F^{p+1}}F^{kl}\br{\bar\n\ov{\Rm}(\nu,x_{;k},x_{;r},x_{;l},x_{;j})+\bar\n\ov{\Rm}(\nu,x_{;r},x_{;j},x_{;k},x_{;l})}g^{ri}.}
\end{lemma}

\pf{In \cref{Ev-W-2} we rewrite the terms involving the Riemann tensor employing \eqref{Warped-Rm}:
\eq{&\bar{R}_{\al\be\g\de}x^{\al}_{;l}x^{\be}_{;j}x^{\g}_{;k}x^{\de}_{;m}\\
	=~&-\fr{\vt''}{\vt}(g_{lm}g_{jk}-g_{lk}g_{jm})+\br{\fr{\vt''}{\vt}-\fr{\vt'^{2}}{\vt^{2}}}(\bar{g}_{lm}\bar{g}_{jk}-\bar{g}_{lk}\bar{g}_{jm})+\ti{R}_{\al\be\g\de}x^{\al}_{;l}x^{\be}_{;j}x^{\g}_{;k}x^{\de}_{;m}\\
			=~&-\fr{\vt'^{2}}{\vt^{2}}(g_{lm}g_{jk}-g_{lk}g_{jm})+\ti{R}_{\al\be\g\de}x^{\al}_{;l}x^{\be}_{;j}x^{\g}_{;k}x^{\de}_{;m}\\
			-&\br{\fr{\vt''}{\vt}-\fr{\vt'^{2}}{\vt^{2}}}(u_{;l}u_{;m}g_{jk}+u_{;j}u_{;k}g_{lm}-u_{;l}u_{;k}g_{jm}-u_{;j}u_{;m}g_{lk}).}
Hence
\eq{&\cF^{kl}\bar{R}_{\al\be\g\de}\br{x^{\al}_{;l}x^{\be}_{;j}x^{\g}_{;k}x^{\de}_{;m}h^{im}+x^{\al}_{;l}x^{\be}_{;r}x^{\g}_{;k}x^{\de}_{;m}h^m_jg^{ri}}\\
	=~&-\fr{\vt'^{2}}{\vt^{2}}\fr{p}{F^{p+1}}F^{kl}h^{i}_{l}g_{jk}+\fr{\vt'^{2}}{\vt^{2}}\fr{p}{F^{p+1}}F^{k}_{k}h^{i}_{j}-\fr{\vt'^{2}}{\vt^{2}}\fr{p}{F^{p+1}}F^{il}h_{jl}+\fr{\vt'^{2}}{\vt^{2}}\fr{p}{F^{p+1}}F^{k}_{k}h^{i}_{j}\\
	+&\br{\fr{\vt''}{\vt}-\fr{\vt'^{2}}{\vt^{2}}}\fr{p}{F^{p+1}}F^{kl}(u_{;k}u_{;l}h^{i}_{j}+g_{kl}h^{im}u_{;m}u_{;j}-h^{i}_{l}u_{;j}u_{;k}-g_{kj}h^{im}u_{;l}u_{;m})\\
	+&\br{\fr{\vt''}{\vt}-\fr{\vt'^{2}}{\vt^{2}}}\fr{p}{F^{p+1}}(F^{kl}(u_{;k}u_{;l}h^{i}_{j}+g_{kl}h^{m}_{j}u_{;m}{u_{;}}^{i})-F^{il}h^{m}_{j}u_{;l}u_{;m}-F^{k}_{m}h^{m}_{j}{u_{;}}^{i}u_{;k})\\
	+&~\fr{p}{F^{p+1}}F^{kl}\br{\ti{R}_{\al\be\g\de}x^{\al}_{;l}x^{\be}_{;j}x^{\g}_{;k}x^{\de}_{;m}h^{im}+\ti{R}_{\al\be\g\de}x^{\al}_{;l}x^{\be}_{;r}x^{\g}_{;k}x^{\de}_{;m}h^{m}_{j}g^{ri}}}
and thus the following two equations hold:
\eq{	&\cF^{kl}\bar{R}_{\al\be\g\de}\br{x^{\al}_{;l}x^{\be}_{;j}x^{\g}_{;k}x^{\de}_{;m}h^{im}+x^{\al}_{;l}x^{\be}_{;r}x^{\g}_{;k}x^{\de}_{;m}h^m_jg^{ri}}\\
	=~&-2\fr{\vt'^{2}}{\vt^{2}}\fr{p}{F^{p+1}}(F^{i}_{l}h^{l}_{j}-F^{kl}g_{kl}h^{i}_{j})\\
	+&\br{\fr{\vt''}{\vt}-\fr{\vt'^{2}}{\vt^{2}}}\fr{p}{F^{p+1}}(2F^{kl}u_{;k}u_{;l}h^{i}_{j}+F^{kl}g_{kl}(h^{im}u_{;m}u_{;j}+h^{m}_{j}u_{;m}{u_{;}}^{i})\\
	&\hp{===========}-F^{l}_{j}h^{im}u_{;l}u_{;m}-F^{il}h^{m}_{j}u_{;l}u_{;m}-F^{kl}h^{i}_{l}u_{;j}u_{;k}-F^{k}_{m}h^{m}_{j}{u_{;}}^{i}u_{;k})\\
	+&~\fr{p}{F^{p+1}}F^{kl}\br{\ti{R}_{\al\be\g\de}x^{\al}_{;l}x^{\be}_{;j}x^{\g}_{;k}x^{\de}_{;m}h^{im}+\ti{R}_{\al\be\g\de}x^{\al}_{;l}x^{\be}_{;r}x^{\g}_{;k}x^{\de}_{;m}h^{m}_{j}g^{ri}},}
\eq{&2\cF^{kl}\bar{R}_{\al\be\g\de}x^{\al}_{;r}x^{\be}_{;m}x^{\g}_{;k}x^{\de}_{;j}h^m_lg^{ri}\\
	=~&-2\fr{\vt'^{2}}{\vt^{2}}\fr{p}{F^{p+1}}(F\de^{i}_{j}-F^{i}_{l}h^{l}_{j})+\fr{2p}{F^{p+1}}F^{kl}\ti{R}_{\al\be\g\de}x^{\al}_{;r}x^{\be}_{;m}x^{\g}_{;k}x^{\de}_{;j}h^{m}_{l}g^{ri}\\
	+&~2\br{\fr{\vt''}{\vt}-\fr{\vt'^{2}}{\vt^{2}}}\fr{p}{F^{p+1}}(F^{k}_{l}h^{l}_{j}{u_{;}}^{i}u_{;k}+F^{il}h^{m}_{l}u_{;j}u_{;m}-F{u_{;}}^{i}u_{;j}-F^{kl}h^{m}_{l}u_{;m}u_{;k}\de^{i}_{j}).}
Adding up, also using $F^{i}_{k}h^{k}_{j}=h^{i}_{k}F^{k}_{j}$, gives
\eq{\label{ICF-Ev-W-1}&\cF^{kl}\bar{R}_{\al\be\g\de}\br{x^{\al}_{;l}x^{\be}_{;j}x^{\g}_{;k}x^{\de}_{;m}h^{im}+x^{\al}_{;l}x^{\be}_{;r}x^{\g}_{;k}x^{\de}_{;m}h^m_jg^{ri}}\\
	+&2\cF^{kl}\bar{R}_{\al\be\g\de}x^{\al}_{;r}x^{\be}_{;m}x^{\g}_{;k}x^{\de}_{;j}h^m_lg^{ri}\\
	=~&2\fr{\vt'^{2}}{\vt^{2}}\fr{p}{F^{p+1}}(F^{k}_{k}h^{i}_{j}-F\de^{i}_{j})\\
	+&\br{\fr{\vt''}{\vt}-\fr{\vt'^{2}}{\vt^{2}}}\fr{p}{F^{p+1}}(2F^{kl}u_{;k}u_{;l}h^{i}_{j}+F^{k}_{k}(h^{i}_{m}{u_{;}}^{m}u_{;j}+h^{m}_{j}u_{;m}{u_{;}}^{i})\\
	&~\hp{==========}-F^{l}_{j}h^{im}u_{;l}u_{;m}-F^{il}h^{m}_{j}u_{;l}u_{;m}+F^{k}_{l}h^{l}_{j}{u_{;}}^{i}u_{;k}\\
	&~\hp{==========}+F^{il}h^{m}_{l}u_{;j}u_{;m}-2F{u_{;}}^{i}u_{;j}-2F^{kl}h^{m}_{l}u_{;m}u_{;k}\de^{i}_{j})\\
	+&~\fr{p}{F^{p+1}}F^{kl}\ti{R}_{\al\be\g\de}\br{x^{\al}_{;l}x^{\be}_{;j}x^{\g}_{;k}x^{\de}_{;m}h^{im}+x^{\al}_{;l}x^{\be}_{;r}x^{\g}_{;k}x^{\de}_{;m}h^{m}_{j}g^{ri}}\\
	+&~\fr{2p}{F^{p+1}}F^{kl}\ti{R}_{\al\be\g\de}x^{\al}_{;r}x^{\be}_{;m}x^{\g}_{;k}x^{\de}_{;j}h^{m}_{l}g^{ri}.}
Using $\nu=v^{-1}(1,-\bar{g}^{ik}u_{;k})$ and $v^{-2}\bar{g}^{ij}u_{;i}u_{;j}=\|\n u\|^{2},$ we get
\eq{\label{ICF-Ev-W-4}\ov{\Rm}(x_{;i},\nu,\nu,x_{;j})&=-\fr{\vt''}{\vt}g_{ij}+\br{\fr{\vt''}{\vt}-\fr{\vt'^{2}}{\vt^{2}}}\br{\|\n u\|^{2}\bar g_{ij}-v^{-2}u_{;i}u_{;j}}\\
								&\hp{=}+\widetilde{\Rm}(x_{;i},\nu,\nu,x_{;j})\\
										&=-\fr{\vt''}{\vt}g_{ij}+\br{\fr{\vt''}{\vt}-\fr{\vt'^{2}}{\vt^{2}}}\br{\|\n u\|^{2}g_{ij}-u_{;i}u_{;j}}\\
						&\hp{=}+\widetilde{\Rm}(x_{;i},\nu,\nu,x_{;j}).}
Thus

\eq{\label{ICF-Ev-W-2}&~(\cF-\cF^{kl}h_{kl})\bar{R}_{\al\be\g\de}x^{\al}_{;r}\nu^{\be}\nu^{\g}x^{\de}_{j}g^{ri}+\cF^{kl}\bar{R}_{\al\be\g\de}x^{\al}_{;k}\nu^{\be}\nu^{\g}x^{\de}_{;l}h^{i}_{j}\\
	=~&\fr{p+1}{F^{p}}\br{\fr{\vt''}{\vt}\de^{i}_{j}-\br{\fr{\vt''}{\vt}-\fr{\vt'^{2}}{\vt^{2}}}(\|\n u\|^{2}\de^{i}_{j}-u_{;j}{u_{;}}^{i})}\\
	-&~\fr{p}{F^{p+1}}F^{kl}\br{\fr{\vt''}{\vt}g_{kl}-\br{\fr{\vt''}{\vt}-\fr{\vt'^{2}}{\vt^{2}}}(\|\n u\|^{2}g_{kl}-u_{;l}{u_{;k}})}h^{i}_{j}\\
	-&~\fr{p+1}{F^{p}}\widetilde{\Rm}(x_{;r},\nu,\nu,x_{;j})g^{ri}+\fr{p}{F^{p+1}}F^{kl}\widetilde{\Rm}(x_{;k},\nu,\nu,x_{;l})h^{i}_{j}.}

%
Adding up \eqref{ICF-Ev-W-1} and \eqref{ICF-Ev-W-2} and inserting the result into \cref{Ev-W-2} gives the claimed formula.
}

\begin{lemma}\label{ICF-Ev-v}
Along \eqref{ICF} the function $v=\vt s^{-1}$
 satisfies the evolution equation
\eq{\label{ICF-Ev-v-1}&\dot{v}-\fr{p}{F^{p+1}}F^{ij}v_{;ij}\\
=~&-\fr{p}{F^{p+1}}F^{ij}h_{ik}h^{k}_{j}v-\fr{\vt'^{2}}{\vt^{2}}\fr{p}{F^{p+1}}F^{ij}g_{ij}v+\fr{\vt'}{\vt}\fr{p+1}{F^{p}}+\fr{\vt'}{\vt}\fr{p-1}{F^{p}}v^{2}\\
					+&~\br{\fr{\vt'^{2}}{\vt^{2}}-\fr{\vt''}{\vt}}\fr{p}{F^{p+1}}F^{ij}u_{;i}u_{;j}v\\
					+&~\fr{p}{F^{p+1}}F^{kl}\br{\fr{\vt''}{\vt}-\fr{\vt'^{2}}{\vt^{2}}}(u_{;l}u_{;k}-\|\n u\|^{2}g_{kl})v\\
				-&~\fr{p}{F^{p+1}}F^{kl}\widetilde{\Rm}(x_{;k},\nu,\nu,x_{;l})v+2\fr{\vt'}{\vt}\fr{p}{F^{p+1}}F^{ij}u_{;i}v_{;j}	-\fr{2}{v}\fr{p}{F^{p+1}}F^{ij}v_{;i}v_{;j}.}
\end{lemma}

\pf{Due to \eqref{Ev-u} and \eqref{ICF-Ev-s}
we have
\eq{&\dot{v}-\fr{p}{F^{p+1}}F^{ij}v_{;ij}\\
	=~&\fr{\vt'}{\vt}v\br{\dot{u}-\fr{p}{F^{p+1}}F^{ij}u_{;ij}}-\fr{\vt''}{\vt}v\fr{p}{F^{p+1}}F^{ij}u_{;i}u_{;j}-\fr{v}{s}\br{\dot{s}-\fr{p}{F^{p+1}}F^{ij}s_{;ij}}\\
					-&~\fr{2\vt}{s^{3}}\fr{p}{F^{p+1}}F^{ij}s_{;i}s_{;j}	-2\fr{p}{F^{p+1}}F^{ij}\vt_{;i}\br{\fr{1}{s}}_{;j}\\
					=~&\fr{\vt'}{\vt}\fr{p+1}{F^{p}}-\fr{\vt'^{2}}{\vt^{2}}\fr{p}{F^{p+1}}F^{ij}g_{ij}v+\br{\fr{\vt'^{2}}{\vt^{2}}-\fr{\vt''}{\vt}}\fr{p}{F^{p+1}}F^{ij}u_{;i}u_{;j}v\\
					-&~\fr{p}{F^{p+1}}F^{ij}h_{ik}h^{k}_{j}v+\fr{\vt'}{\vt}\fr{p-1}{F^{p}}v^{2}-\fr{p}{F^{p+1}}F^{kl}\ov{\Rm}(\nu,x_{;k},x_{;m},x_{;l}){u_{;}}^{m}v^{2}\\
					-&~\fr{2}{v}\fr{p}{F^{p+1}}F^{ij}v_{;i}v_{;j}	-2v\fr{\vt'^{2}}{\vt^{2}}\fr{p}{F^{p+1}}F^{ij}u_{;i}u_{;j}+4\fr{\vt'}{\vt}\fr{p}{F^{p+1}}F^{ij}u_{;i}v_{;j}\\
					-&~2\fr{\vt'}{\vt}\fr{p}{F^{p+1}}F^{ij}u_{;i}v_{;j}+2\fr{\vt'^{2}}{\vt^{2}}\fr{p}{F^{p+1}}F^{ij}u_{;i}u_{;j}v,}
which is the claimed formula up to rewriting the term involving $\ov{\Rm}$. However, we use \eqref{Warped-Rm} to deduce
\eq{\ov{\Rm}(\nu,x_{;k},x_{;m},x_{;l}){u_{;}}^{m}&=-\br{\fr{\vt''}{\vt}-\fr{\vt'^{2}}{\vt^{2}}}v^{-1}(u_{;l}u_{;k}-\|\n u\|^{2}g_{kl})\\
			&\hp{=}+\widetilde{\Rm}(\nu,x_{;k},x_{;m},x_{;l}){u_{;}}^{m}\\
			}
and hence
\eq{\ov{\Rm}(\nu,x_{;k},x_{;m},x_{;l}){u_{;}}^{m}v^{2}&=-\br{\fr{\vt''}{\vt}-\fr{\vt'^{2}}{\vt^{2}}}(u_{;l}u_{;k}-\|\n u\|^{2}g_{kl})v\\
				&\hp{=}+\widetilde{\Rm}(x_{;k},\nu,\nu,x_{;l})v,}
where we have used
\eq{{u_{;}}^{m}=v^{-2}\vt^{-2}\s^{mk}u_{;k}.}
Inserting gives the result.
}

We start the investigation of the long-time behavior of \eqref{ICF-Ini} under the assumptions in item (ii) of \cref{ICF-Main} by proving a lower bound on the curvature function.

\begin{lemma}\label{ICF-LT-1/F-Bound}
Under the assumptions of \cref{ICF-Main}~(ii), along \eqref{ICF} there exists a constant $c$, such that
\eq{\fr{\vt'v}{\vt F}\leq c.}
\end{lemma}

\pf{
If $p=1$ and $\vt'$ is bounded , the result follows from \cref{ICF-1/F-bound} immediately. If $p<1$ or $\vt'$ is unbounded, \cref{ICF-grad-bound} says that $v\ra 1$.
Due \eqref{Ev-F}, \eqref{Ev-u}, \eqref{ICF-Ev-W-4} and \eqref{ICF-Ev-v-1}  the function
\eq{w=\log\br{\fr{1}{F^{p}}}+f(v)+p\log \vt'-p\log\vt,}
where $f$ with $f'\geq 0$ is yet to be determined, satisfies
\eq{&\dot{w}-\fr{p}{F^{p+1}}F^{ij}w_{;ij}\\
	=~&\fr{p}{F^{p+1}}F^{ij}h_{ik}h^{k}_{j}-\fr{\vt''}{\vt}\fr{p}{F^{p+1}}F^{ij}g_{ij}+\fr{p}{F^{p+1}}F^{ij}\widetilde{\Rm}(x_{;i},\nu,\nu,x_{;j})\\
	+&\br{\fr{\vt''}{\vt}-\fr{\vt'^{2}}{\vt^{2}}}\fr{p}{F^{p+1}}F^{ij}(\|\n u\|^{2}g_{ij}-u_{;i}u_{;j})\\
	+&~\fr{p}{F^{p+1}}F^{ij}\br{\log\fr{1}{F^{p}}}_{;i}\br{\log\fr{1}{F^{p}}}_{;j}-\fr{p}{F^{p+1}}F^{ij}h_{ik}h^{k}_{j}f'v-\fr{\vt'^{2}}{\vt^{2}}\fr{p}{F^{p+1}}F^{ij}g_{ij}f'v\\
	+&~\fr{\vt'}{\vt}\fr{p+1}{F^{p}}f'+\fr{\vt'}{\vt}\fr{p-1}{F^{p}}f'v^{2}+\br{\fr{\vt'^{2}}{\vt^{2}}-\fr{\vt''}{\vt}}\fr{p}{F^{p+1}}F^{ij}u_{;i}u_{;j}f'v\\
	+&\br{\fr{\vt''}{\vt}-\fr{\vt'^{2}}{\vt^{2}}}\fr{p}{F^{p+1}}F^{ij}(u_{;i}u_{;j}-\|\n u\|^{2}g_{ij})f'v\\
	-&~\fr{p}{F^{p+1}}F^{ij}\widetilde{\Rm}(x_{;i},\nu,\nu,x_{;j})f'v+2f'\fr{\vt'}{\vt}\fr{p}{F^{p+1}}F^{ij}u_{;i}v_{;j}-\fr{2}{v}f'\fr{p}{F^{p+1}}F^{ij}v_{;i}v_{;j}\\
	-&~f''\fr{p}{F^{p+1}}F^{ij}v_{;i}v_{;j}+(p+1)\br{\fr{\vt''}{\vt'}-\fr{\vt'}{\vt}}\fr{p}{F^{p+1}}v^{-1}F\\
	-&~p\br{\fr{\vt''}{\vt}-\fr{\vt'^{2}}{\vt^{2}}}\fr{p}{F^{p+1}}F^{ij}g_{ij}+p\br{\fr{\vt''}{\vt}-\fr{\vt'^{2}}{\vt^{2}}}\fr{p}{F^{p+1}}F^{ij}u_{;i}u_{;j}\\
	-&~p\br{\fr{\vt''}{\vt'}-\fr{\vt'}{\vt}}'\fr{p}{F^{p+1}}F^{ij}u_{;i}u_{;j}.}

Sorting the terms appropriately and replacing $(\log 1/F^{p})_{;i}$ we get
\eq{\label{ICF-LT-1/F-Bound-1}&\dot{w}-\fr{p}{F^{p+1}}F^{ij}w_{;ij}\\
	\leq~&\fr{p}{F^{p+1}}F^{ij}h_{ik}h^{k}_{j}(1-f'v)+\fr{p}{F^{p+1}}F^{ij}\widetilde{\Rm}(x_{;i},\nu,\nu,x_{;j})(1-f'v)\\
	+&\br{\fr{\vt''}{\vt}-\fr{\vt'^{2}}{\vt^{2}}}\fr{p}{F^{p+1}}F^{ij}(\|\n u\|^{2}g_{ij}-u_{;i}u_{;j})(1-f'v)\\
	+&~\fr{p(p+1)}{F^{p+1}}\fr{\vt''}{\vt}\br{-F^{ij}g_{ij}+\fr{\vt}{\vt'}Fv^{-1}+c(f'v+1)F^{ij}u_{;i}u_{;j}}\\
	+&~\fr{\vt'^{2}}{\vt^{2}}\fr{p}{F^{p+1}}\Big(-F^{ij}g_{ij}f'v+\fr{\vt}{\vt'}\fr{p+1}{p}f'F+\fr{\vt}{\vt'}\fr{p-1}{p}f'Fv^{2}\\
	-&~(p+1)\fr{\vt}{\vt'}Fv^{-1}	+pF^{ij}g_{ij}\Big)+\fr{\vt'^{2}}{\vt^{2}}\fr{p}{F^{p+1}}F^{ij}u_{;i}u_{;j}(f'v-1)\\
	+&~(p-1)^{2}\fr{\vt'^{2}}{\vt^{2}}\fr{p}{F^{p+1}}F^{ij}u_{;i}u_{;j}+\fr{p}{F^{p+1}}F^{ij}w_{;i}w_{;j}-2f'\fr{p}{F^{p+1}}F^{ij}w_{;i}v_{;j}\\
	-&~2p\br{\fr{\vt''}{\vt'}-\fr{\vt'}{\vt}}\fr{p}{F^{p+1}}F^{ij}w_{;i}u_{;j}+f'^{2}\fr{p}{F^{p+1}}F^{ij}v_{;i}v_{;j}\\
	-&~2pf'\br{\fr{\vt''}{\vt'}-\fr{\vt'}{\vt}}\fr{p}{F^{p+1}}F^{ij}v_{;i}u_{;j}+2f'\fr{\vt'}{\vt}\fr{p}{F^{p+1}}F^{ij}u_{;i}v_{;j}\\
	-&~\fr{2}{v}f'\fr{p}{F^{p+1}}F^{ij}v_{;i}v_{;j}-f''\fr{p}{F^{p+1}}F^{ij}v_{;i}v_{;j}.}
Now choose 
\eq{f(v)=-\log\br{v^{-\fr 12}-\fr 34}}
where in the sequel we only consider sufficiently large times where $v^{\fr12}<\fr 43$. Then
\eq{f'=\fr{1}{2}\fr{v^{-\fr 32}}{v^{-\fr 12}-\fr 34},\q f''=-\fr 34\fr{v^{-\fr 52}}{v^{-\fr 12}-\fr 34}+\fr 14\fr{v^{-3}}{\br{v^{-\fr 12}-\fr 34}^{2}},\q f'v=\fr 12\fr{v^{-\fr 12}}{v^{-\fr 12}-\fr 34}\geq \fr 32}
and 
\eq{f'^{2}-\fr{2}{v}f'-f''&=-\fr{v^{-\fr 52}}{v^{-\fr 12}-\fr 34}+\fr 34\fr{v^{-\fr 52}}{v^{-\fr 12}-\fr 34}\\
				&=-\fr{1}{4}\fr{v^{-\fr 52}}{v^{-\fr 12}-\fr 34}\\
				&\leq -\fr{3}{4v^{2}}.}
				
Hence, using Cauchy-Schwarz on $F^{ij}u_{;i}v_{;j}$, we can estimate further:

\eq{&\dot{w}-\fr{p}{F^{p+1}}F^{ij}w_{;ij}\\
	\leq~&\fr{p(p+1)}{F^{p+1}}\fr{\vt''}{\vt}\br{-F^{ij}g_{ij}+\fr{\vt}{\vt'}Fv^{-1}+c_{\e}(f'v+1)F^{ij}u_{;i}u_{;j}}\\
	+&~\fr{\vt'^{2}}{\vt^{2}}\fr{p}{F^{p+1}}\Big(pF^{ij}g_{ij}-F^{ij}g_{ij}f'v+2\fr{\vt}{\vt'}f'F-(p+1)\fr{\vt}{\vt'}Fv^{-1}+F^{ij}g_{ij}\|\n u\|^{2}f'v\\
	+&~c_{\e}(f'v+1)F^{ij}u_{;i}u_{;j}\Big)+\fr{p}{F^{p+1}}F^{ij}v_{;i}v_{;j}\br{-\fr{3}{4v^{2}}+\fr{\e c}{v}f'}\\
	+&~\fr{p}{F^{p+1}}F^{ij}w_{;i}w_{;j}-2f'\fr{p}{F^{p+1}}F^{ij}w_{;i}v_{;j}-2p\br{\fr{\vt''}{\vt'}-\fr{\vt'}{\vt}}\fr{p}{F^{p+1}}F^{ij}w_{;i}u_{;j}\\
	<~&0}
at maximal points of $w$, if $\e$ is chosen small, $\|\n u\|$ is small enough (which happens eventually) and $w$ is large. Hence $w$ is bounded.
}

We need a similar estimate of the rescaled principal curvatures.

\begin{lemma}\label{ICF-Scaled-W}
Under the assumptions of \cref{ICF-Main}~(ii), along \eqref{ICF} there exists a constant $c$, such that
\eq{\ka_{n}\leq c\fr{\vt'}{\vt}.}
\end{lemma}

\pf{
Define
\eq{z=\log h^{n}_{n}+\log\fr{\vt}{\vt'}+f(v).}
 We estimate the evolution of $z$ directly from \eqref{Ev-u}, \cref{ICF-Ev-W} and \eqref{ICF-Ev-v-1} and, as in the proof of \cref{ICF-k-bound}, from the start calculate in a maximal point of $z$ in coordinates such that
\eq{g_{ij}=\de_{ij},\q h_{ij}=\ka_{i}\de_{ij},\q \ka_{1}\leq \dots\leq \ka_{n}.}
First of all there holds
\eq{&F^{kl}\wt{\Rm}(x_{;l},x_{;j},x_{;k},x_{;m})h^{im}+F^{kl}\wt{\Rm}(x_{;l},x_{;r},x_{;k},x_{;m})h^{m}_{j}g^{ri}\\
	+~&2F^{kl}\wt{\Rm}(x_{;r},x_{;m},x_{;k},x_{;j})h^{m}_{l}g^{ri}\\
	=~& 2F^{kk}\wt{\Rm}(x_{;n},x_{;k},x_{;k},x_{;n})(\ka_{k}-\ka_{n})\\
	\leq~&0,}
since the sectional curvatures of $\s$ are non-negative.

Due to \cref{tildeRm} we have
\eq{\|\widetilde{\Rm}\|\leq \fr{c}{\vt^{2}},\q \|\bar\n\widetilde{\Rm}\|\leq c\fr{\vt'}{\vt^{3}}} 
and we get

\eq{&F^{kl}\bar{\n}\ov{\Rm}(\nu,x_{;k},x_{;r},x_{;l},x_{;j})g^{ri}+F^{kl}\bar{\n}\ov{\Rm}(\nu,x_{;r},x_{;j},x_{;k},x_{;l})g^{ri}\\
	\leq~&c\fr{\vt'^{3}}{\vt^{3}}\|\n u\|^{2}F^{k}_{k}+c\fr{\vt'}{\vt^{3}}\|\n u\|F^{k}_{k},}
where we have used that the terms in \eqref{Warped-Rm-Der-1} involving $r_{;\al}$ are cancelled, since $\bar{T}$ carries the symmetries of a curvature tensor.	

Hence

\eq{&\dot{z}-\fr{p}{F^{p+1}}F^{ij}z_{;ij}\\
	\leq ~&\fr{p}{F^{p+1}}F^{ij}(\log h^{n}_{n})_{;i}(\log h^{n}_{n})_{;j}+\fr{p}{F^{p+1}}F^{ij}h_{ik}h^{k}_{j}(1-f'v)-\fr{p+1}{F^{p}}h^{n}_{n}\\
	+&~\fr{\vt'^{2}}{\vt^{2}}\fr{p}{F^{p+1}}F^{k}_{k}(1-f'v)-\fr{\vt'^{2}}{\vt^{2}}\fr{p-1}{F^{p}}\ka_{n}^{-1}-\br{\fr{\vt''}{\vt}-\fr{\vt'^{2}}{\vt^{2}}}\fr{p}{v^{2}F^{p+1}}F^{k}_{k}\\
	+&\br{\fr{\vt''}{\vt}-\fr{\vt'^{2}}{\vt^{2}}}\fr{p+1}{v^{2}F^{p}}\ka_{n}^{-1}+c(1+|f'|)\fr{\vt'^{2}}{\vt^{2}}\fr{p}{F^{p+1}}F^{k}_{k}\|\n u\|^{2}\\
	+&~c\fr{\vt'^{2}}{\vt^{2}}\fr{1}{F^{p}}\|\n u\|^{2}\ka_{n}^{-1}+\fr{c}{F^{p+1}}\fr{\vt'^{3}}{\vt^{3}}\|\n u\|^{2}\ka_{n}^{-1}F^{k}_{k}+\fr{c}{F^{p+1}}\fr{\vt'}{\vt^{3}}\|\n u\|\ka_{n}^{-1}F^{k}_{k}\\
	+&~f'\fr{\vt'}{\vt}\fr{p+1}{F^{p}}+f'\fr{\vt'}{\vt}\fr{p-1}{F^{p}}v^{2}+2f'\fr{\vt'}{\vt}\fr{p}{F^{p+1}}F^{ij}u_{;i}v_{;j}-\fr{2f'}{v}\fr{p}{F^{p+1}}F^{ij}v_{;i}v_{;j}\\
	-&~f''\fr{p}{F^{p+1}}F^{ij}v_{;i}v_{;j}+\br{\fr{\vt'}{\vt}-\fr{\vt''}{\vt'}}\fr{p+1}{F^{p}}v^{-1}-\br{\fr{\vt'^{2}}{\vt^{2}}-\fr{\vt''}{\vt}}\fr{p}{F^{p+1}}F^{k}_{k}.}
Hence

\eq{\label{ICF-Scaled-W-3}&\dot{z}-\fr{p}{F^{p+1}}F^{ij}z_{;ij}\\
	\leq ~&\fr{p}{F^{p+1}}F^{ij}(\log h^{n}_{n})_{;i}(\log h^{n}_{n})_{;j}+\fr{p}{F^{p+1}}F^{ij}h_{ik}h^{k}_{j}(1-f'v)\\
	+&~\fr{\vt'^{2}}{\vt^{2}}\fr{p}{F^{p+1}}F^{k}_{k}(1-f'v)-\fr{p+1}{F^{p}}\br{\ka_{n}-\fr{\vt'}{\vt}f'}\\
	+&\br{\fr{\vt''}{\vt'}-\fr{\vt'}{\vt}}\fr{p+1}{F^{p}}v^{-1}\br{v^{-1}\ka_{n}^{-1}\fr{\vt'}{\vt}-1}+\fr{\vt'}{\vt}\fr{p-1}{F^{p}}\br{f'v^{2}-\fr{\vt'}{\vt}\ka_{n}^{-1}}\\
	+&~\fr{\vt'^{2}}{\vt^{2}}\br{\fr{\vt'}{\vt}\ka_{n}^{-1}+1}\fr{c(1+|f'|)}{F^{p+1}}F^{k}_{k}\|\n u\|^{2}+\fr{\vt'^{2}}{\vt^{2}}\fr{c}{F^{p}}\|\n u\|^{2}\ka_{n}^{-1}\\
	+&~\fr{c}{F^{p+1}}\fr{\vt'}{\vt^{3}}\|\n u\|\ka_{n}^{-1}F^{k}_{k}+2f'\fr{\vt'}{\vt}\fr{p}{F^{p+1}}F^{ij}u_{;i}v_{;j}-\fr{2f'}{v}\fr{p}{F^{p+1}}F^{ij}v_{;i}v_{;j}\\
	-&~f''\fr{p}{F^{p+1}}F^{ij}v_{;i}v_{;j}. }
Furthermore we insert 
\eq{(\log h^{n}_{n})_{;i}=z_{;i}-\br{\fr{\vt'}{\vt}-\fr{\vt''}{\vt'}}u_{;i}-f'v_{;i}}
and hence
\eq{\label{ICF-Scaled-W-2}&\dot{z}-\fr{p}{F^{p+1}}F^{ij}z_{;ij}\\
	\leq ~&\fr{p}{F^{p+1}}F^{ij}h_{ik}h^{k}_{j}(1-f'v)+\fr{\vt'^{2}}{\vt^{2}}\fr{p}{F^{p+1}}F^{k}_{k}(1-f'v)-\fr{p+1}{F^{p}}\br{\ka_{n}-\fr{\vt'}{\vt}f'}\\
	+&~\br{\fr{\vt''}{\vt'}-\fr{\vt'}{\vt}}\fr{p+1}{F^{p}}v^{-1}\br{v^{-1}\ka_{n}^{-1}\fr{\vt'}{\vt}-1}+\fr{p-1}{F^{p}}\br{\fr{\vt'}{\vt}f'v^{2}-\fr{\vt'^{2}}{\vt^{2}}\ka_{n}^{-1}}\\
	+&~\fr{\vt'^{2}}{\vt^{2}}\br{\fr{\vt'}{\vt}\ka_{n}^{-1}+1}\fr{c(1+|f'|)}{F^{p+1}}F^{k}_{k}\|\n u\|^{2}+\fr{\vt'^{2}}{\vt^{2}}\fr{c}{F^{p}}\|\n u\|^{2}\ka_{n}^{-1}\\
	+&~\fr{c}{F^{p+1}}\fr{\vt'}{\vt^{3}}\|\n u\|\ka_{n}^{-1}F^{k}_{k}+2f'\fr{\vt'}{\vt}\fr{p}{F^{p+1}}F^{ij}u_{;i}v_{;j}\\
	+&~\fr{p}{F^{p+1}}F^{ij}v_{;i}v_{;j}\br{f'^{2}-\fr{2f'}{v}-f''}+2\br{\fr{\vt'}{\vt}-\fr{\vt''}{\vt'}}f'\fr{p}{F^{p+1}}F^{ij}u_{;i}v_{;j}\\
	+&~\fr{p}{F^{p+1}}F^{ij}z_{;i}z_{;j}-2\fr{p}{F^{p+1}}F^{ij}z_{;i}\br{f'v_{;j}+\br{\fr{\vt'}{\vt}-\fr{\vt''}{\vt'}}u_{;j}}. }
Pick
\eq{f(v)=-\log(v^{-\al}-\be),}
where
 \eq{0<\be<\fr{1}{2v},\q 1-\fr{\be}{2}<\al<1. }
 Then
 \eq{1-f'v=\fr{(1-\al)v^{-\al}-\be}{v^{-\al}-\be}\leq \fr{\be(\fr{v^{-\al}}{2}-1)}{v^{-\al}-\be}\leq -\fr{\be}{2}<0}
 and 
 \eq{f'^{2}-\fr{2}{v}f'-f''=\fr{\al v^{-(\al+2)}}{v^{-\al}-\be}(\al-1)\leq\fr{3}{4}\fr{\al-1}{v^{2}}<0.}
  Hence at a maximal point of $z$ there holds
\eq{\label{ICF-Scaled-W-1}&\dot{z}-\fr{p}{F^{p+1}}F^{ij}z_{;ij}\\
	\leq ~&\fr{1}{F^{p}}\br{-(p+1)h^{n}_{n}+c\fr{\vt'}{\vt}+c\fr{\vt'^{2}}{\vt^{2}}(h^{n}_{n})^{-1}}\\
	+&~\fr{\vt'^{2}}{\vt^{2}}\fr{p}{F^{p+1}}F^{k}_{k}\br{-\fr{\be}{2}+c_{\e}\br{1+\fr{\vt'}{\vt}\ka_{n}^{-1}}\|\n u\|^{2}+c\fr{\vt'}{\vt}\fr{\|\n u\|}{\vt'^{2}}\ka_{n}^{-1}}\\
	+&~\br{\fr{3}{4}\fr{\al-1}{v^{2}}+\e c f'}\fr{p}{F^{p+1}}F^{ij}v_{;i}v_{;j},}
where we used
\eq{2|F^{ij}u_{;i}v_{;j}|\leq \fr{\e\vt}{\vt'}F^{ij}v_{;i}v_{;j}+\fr{\vt'}{\e\vt}F^{ij}u_{;i}u_{;j}}
with sufficiently small $\e$.
In case $\sup_{r>0}\vt'(r)=\8$ or $p<1$, we have $\|\n u\|^{2}\ra 0$ and hence the result follows from the maximum principle. If $\vt'\leq c$ and $p=1$ we supposed that 
\eq{F=n\fr{H_{k+1}}{H_{k}}}
which implies $F^{k}_{k}\leq c$, cf. \cite[Lemma~2.7]{Lu:09/2016}, and hence 
\eq{\fr{\vt'^{2}}{\vt^{2}}\fr{p}{F^{p+1}}F^{k}_{k}\leq c\fr{\vt'}{\vt}\fr{1}{F^{p}}}
due to \cref{ICF-LT-1/F-Bound}. Hence the term $-(p+1)h^{n}_{n}$ dominates the whole evolution and we also obtain the bound on $z$ in this case.
}

\begin{rem}
\cref{ICF-Scaled-W} is the only place where we need that $F$ has this special form in case of bounded $\vt'$. Of course the Euclidean case is excluded from this restriction, since the error terms involving $\|\n u\|^{2}$ will not appear here. However, the Euclidean case has already been settled in \cite{Gerhardt:/1990}.
\end{rem}

\subsection{Decay}

The global bounds from \cref{ICF-LT-1/F-Bound} and \cref{ICF-Scaled-W} as well as 
\eq{F_{|\del\G}=0}
imply that the rescaled principal curvatures
\eq{\ti{\ka}_{i}=\fr{\vt}{\vt'}\ka_{i}}
range in a compact subset of $\G$ and hence the elliptic operator $d_{h}F$ is uniformly bounded,
\eq{c\|\xi\|^{2}\leq d_{h}F(\xi,\xi)\leq C\|\xi\|^{2}.}

The aim of this final section is to show that all $\ti{\ka}_{i}$ actually behave according to the convergence rates described in \cref{ICF-Main}. The following two lemmata prepare this result. Throughout this whole section, the procedure is similar to the one in \cite{Scheuer:01/2017}.

\begin{lemma}\label{ICF-LT-1/F-Decay}
Under the assumptions of \cref{ICF-Main}~(ii), along \eqref{ICF} there exist constants $\mu$ and $c$, such that
\eq{\fr{\vt'}{\vt}{\fr{1}{F}}-\fr{1}{n}\leq \fr{c}{\vt'^{\mu}}}
\end{lemma}

\pf{We only have to consider the case that $\vt'$ is unbounded. Come back to the proof of \cref{ICF-LT-1/F-Bound} and consider \eqref{ICF-LT-1/F-Bound-1} with 
\eq{f(v)=\log v.}

Hence from \eqref{ICF-LT-1/F-Bound-1} we deduce that
\eq{z=\log\br{\fr{1}{F^{p}}}+\log v+p\log \vt'-p\log\vt+p\log n}
satisfies
\eq{&\dot{z}-\fr{p}{F^{p+1}}F^{ij}z_{;ij}\\
\leq ~& \fr{p(p+1)}{F^{p+1}}\fr{\vt''}{\vt}\br{-F^{ij}g_{ij}+\fr{\vt}{\vt'}Fv^{-1}+c\|\n u\|^{2}}\\
	+&~ \fr{\vt'^{2}}{\vt^{2}}\fr{p}{F^{p+1}}\br{(p-1)F^{ij}g_{ij}+\fr{1-p^{2}}{p}\fr{\vt}{\vt'}Fv^{-1}+\fr{\vt}{\vt'}\fr{p-1}{p}Fv+c\|\n u\|^{2}}\\
	+&~\fr{p}{F^{p+1}}F^{ij}z_{;i}z_{;j}-\fr{2}{v}\fr{p}{F^{p+1}}F^{ij}z_{;i}v_{;j}-2p\br{\fr{\vt''}{\vt'}-\fr{\vt'}{\vt}}\fr{p}{F^{p+1}}F^{ij}z_{;i}u_{;j}\\
	\leq~&\fr{np(p+1)}{F^{p+1}}\fr{\vt''}{\vt}\br{-1+e^{-\fr zp}+c\|\n u\|^{2}}+\fr{np(1-p)}{F^{p+1}}\fr{\vt'^{2}}{\vt^{2}}\br{-1+e^{-\fr zp}+c\|\n u\|^{2}}\\
	+&~\fr{p}{F^{p+1}}F^{ij}z_{;i}z_{;j}-\fr{2}{v}\fr{p}{F^{p+1}}F^{ij}z_{;i}v_{;j}-2p\br{\fr{\vt''}{\vt'}-\fr{\vt'}{\vt}}\fr{p}{F^{p+1}}F^{ij}z_{;i}u_{;j}.
	}
For $\mu\geq 0$ define
\eq{\rho=(e^{z}-1)\vt'^{\mu}.}
Then

\eq{&\dot{\rho}-\fr{p}{F^{p+1}}F^{ij}\rho_{;ij}\\
=~&\br{\dot{z}-\fr{p}{F^{p+1}}z_{;ij}}e^{z}\vt'^{\mu}-\fr{p}{F^{p+1}}F^{ij}z_{;i}z_{;j}e^{z}\vt'^{\mu}+\mu\fr{\vt''}{\vt'}\br{\dot{u}-\fr{p}{F^{p+1}}F^{ij}u_{;i}u_{;j}}\rho\\
-&~\br{\mu(\mu-1)\fr{\vt''^{2}}{\vt'^{2}}+\mu\fr{\vt'''}{\vt'}}\fr{p}{F^{p+1}}F^{ij}u_{;i}u_{;j}\rho\\
\leq~&\fr{np(p+1)}{F^{p+1}}\fr{\vt''}{\vt}e^{\fr{p-1}{p}z}\Big((1-e^{\fr zp})\vt'^{\mu}+\fr{\mu}{np}F\fr{\vt}{\vt'}v^{-1}e^{\fr{1-p}{p}z}\rho\\
		&\hp{\fr{np(p+1)}{F^{p+1}}\fr{\vt''}{\vt}e^{\fr{p-1}{p}z}}-\fr{\mu}{n(p+1)}F^{ij}g_{ij}e^{\fr{1-p}{p}z}\rho+c_{\mu}\|\n u\|^{2}\vt'^{\mu}\Big)\\
+&~\fr{np(1-p)}{F^{p+1}}\fr{\vt'^{2}}{\vt^{2}}e^{\fr{p-1}{p}z}\br{-e^{\fr zp}+1+c\|\n u\|^{2}}\vt'^{\mu}-\fr{2}{v}\fr{p}{F^{p+1}}F^{ij}z_{;i}v_{;j}e^{z}\vt'^{\mu}\\
-&~2p\br{\fr{\vt''}{\vt'}-\fr{\vt'}{\vt}}\fr{p}{F^{p+1}}F^{ij}z_{;i}u_{;j}e^{z}\vt'^{\mu}.
}

Now we estimate at maximal points of $\rho$ and thus may assume $z>0$. Then, also using
\eq{0=\vt'^{-\mu}\rho_{;i}&=z_{;i}e^{z}+\mu\fr{\vt''}{\vt'}(e^{z}-1)u_{;i},}
we obtain
\eq{\label{ICF-LT-1/F-Decay-1}&\dot{\rho}-\fr{p}{F^{p+1}}F^{ij}\rho_{;ij}\\
	\leq~&\fr{np(p+1)}{F^{p+1}}\fr{\vt''}{\vt}e^{\fr{p-1}{p}z}\Big(-\rho+\fr{\mu}{np}F\fr{\vt}{\vt'}v^{-1}e^{\fr{1-p}{p}z}\rho\\
		&\hp{\fr{np(p+1)}{F^{p+1}}\fr{\vt''}{\vt}e^{\fr{p-1}{p}z}}-\fr{\mu}		{n(p+1)}F^{ij}g_{ij}e^{\fr{1-p}{p}z}\rho+c\|\n u\|^{2}\vt'^{\mu}\Big)\\
	+&~\fr{np(1-p)}{F^{p+1}}\fr{\vt'^{2}}{\vt^{2}}e^{\fr{p-1}{p}z}\br{-\rho+c\|\n u\|^{2}\vt'^{\mu}},\\
}
which is negative for sufficiently small $\mu$ and large times, due to \cref{ICF-grad-bound} and the remarks at the beginning of this section. The proof is complete.}

\begin{lemma}\label{ICF-LT-k-Decay}
Under the assumptions of \cref{ICF-Main}~(ii), along \eqref{ICF} the $i$-th rescaled principal curvature converges uniformly to $1$,
\eq{\left|v\kappa_{i}\fr{\vt}{\vt'}-1\right|\ra 0,}
provided $\vt'$ is unbounded.
\end{lemma}

\pf{Using \eqref{ICF-Scaled-W-3} with $f(v)=\log v$ we obtain that
\eq{z=\log h^{n}_{n}+\log\fr{\vt}{\vt'}+\log v}
satisfies
\eq{&\dot{z}-\fr{p}{F^{p+1}}F^{ij}z_{;ij}\\
	\leq~&\fr{p}{F^{p+1}}F^{ij}(\log h^{n}_{n})_{;i}(\log h^{n}_{n})_{;j}+c\fr{\vt'^{2}}{\vt^{2}}\fr{1}{F^{p+1}}\|\n u\|^{2}-\fr{p+1}{F^{p}}\fr{\vt'}{\vt}v^{-1}\br{e^{z}-1}\\
	+&~\br{\fr{\vt''}{\vt'}-\fr{\vt'}{\vt}}\fr{p+1}{F^{p}}v^{-1}\br{e^{-z}-1}+\fr{\vt'}{\vt}\fr{p-1}{F^{p}}v\br{1-e^{-z}}+\fr{c}{F^{p+1}}\fr{1}{\vt^{2}}\|\n u\|\\
	=~&\fr{p}{F^{p+1}}F^{ij}(\log h^{n}_{n})_{;i}(\log h^{n}_{n})_{;j}-\fr{\vt''}{\vt'}\fr{p+1}{F^{p}}v^{-1}(1-e^{-z})+c\fr{\vt'^{2}}{\vt^{2}}\fr{1}{F^{p+1}}\|\n u\|^{2}\\
	-&~\fr{\vt'}{\vt}\fr{p+1}{F^{p}}v^{-1}e^{-z}(e^{z}-1)^{2}
	+\fr{\vt'}{\vt}\fr{p-1}{F^{p}}v\br{1-e^{-z}}+\fr{c}{F^{p+1}}\fr{1}{\vt^{2}}\|\n u\|.}
Define
\eq{\rho=(e^{z}-1)\vt'^{\mu},}
with $\mu\geq 0$.
$\rho$ satisfies
\eq{&\dot{\rho}-\fr{p}{F^{p+1}}F^{ij}\rho_{;ij}\\
=~&\br{\dot{z}-\fr{p}{F^{p+1}}z_{;ij}}e^{z}\vt'^{\mu}-\fr{p}{F^{p+1}}F^{ij}z_{;i}z_{;j}e^{z}\vt'^{\mu}+\mu\fr{\vt''}{\vt'}\br{\dot{u}-\fr{p}{F^{p+1}}F^{ij}u_{;i}u_{;j}}\rho\\
-&\br{\mu(\mu-1)\fr{\vt''^{2}}{\vt'^{2}}+\mu\fr{\vt'''}{\vt'}}\fr{p}{F^{p+1}}F^{ij}u_{;i}u_{;j}\rho\\
\leq~& -\fr{\vt''}{\vt'}\fr{p+1}{F^{p}}v^{-1}\rho-\fr{\vt'}{\vt}\fr{p+1}{F^{p}}v^{-1}(e^{z}-1)\rho-\fr{\vt'}{\vt}\fr{1-p}{F^{p}}v\rho+\fr{\vt'^{2}}{\vt^{2}}\fr{c}{F^{p+1}}\|\n u\|^{2}\vt'^{\mu}\\
+&~\fr{c}{F^{p+1}}\fr{1}{\vt^{2}}\|\n u\|\vt'^{\mu}+\fr{p}{F^{p+1}}F^{ij}(\log h^{n}_{n})_{;i}(\log h^{n}_{n})_{;j}e^{z}\vt'^{\mu}\\
-&~\fr{p}{F^{p+1}}F^{ij}z_{;i}z_{;j}e^{z}\vt'^{\mu}+\mu\fr{\vt''}{\vt'}\fr{p+1}{F^{p}}v^{-1}\rho-\mu\fr{\vt''}{\vt}\fr{p}{F^{p+1}}F^{ij}g_{ij}\rho.}
At spatial maxima of $\rho$ we have
\eq{0=\vt'^{-\mu}\rho_{;i}&=z_{;i}e^{z}+\mu\fr{\vt''}{\vt'}(e^{z}-1)u_{;i}\\
			&=h^{n}_{n;i}v\fr{\vt}{\vt'}+h_{n}^{n}\br{\fr{\vt}{\vt'}}_{;i}v+h^{n}_{n}\fr{\vt}{\vt'}v_{;i}+\mu\fr{\vt''}{\vt'}(e^{z}-1)u_{;i}}
and hence
\eq{F^{ij}(\log h^{n}_{n})_{;i}(\log h^{n}_{n})_{;j}\leq c\fr{\vt'^{2}}{\vt^{2}}\|\n u\|^{2}.}
We obtain at a maximal point where $\rho>0$
\eq{\label{ICF-LT-k-Decay-1}&\dot{\rho}-\fr{p}{F^{p+1}}F^{ij}\rho_{;ij}\\
\leq~& -\fr{\vt'}{\vt}\fr{p+1}{F^{p}}v^{-1}(e^{z}-1)\rho-\fr{\vt'}{\vt}\fr{1-p}{F^{p}}v\rho\\
	+&~\fr{n}{F^{p}}\fr{\vt''}{\vt'}v^{-1}\rho\br{-\mu p\fr{v}{F}\fr{\vt'}{\vt}+(p+1)\fr{\mu-1}{ n}}+\fr{\vt'^{2}}{\vt^{2}}\fr{c}{F^{p+1}}\|\n u\|^{2}\vt'^{\mu}\\
	+&~\fr{c}{F^{p+1}}\fr{1}{\vt^{2}}\|\n u\|\vt'^{\mu}\\
	\leq~& \fr{\vt'}{\vt}\fr{p+1}{vF^{p}}\br{-\fr{1-p}{p+1}\rho-(e^{z}-1)\rho+c\|\n u\|^{2}\vt'^{\mu}+\fr{c\vt'^{\mu}}{\vt'^{2}}\|\n u\|}\\
	+&~\fr{ n}{F^{p}}\fr{\vt''}{\vt'}v^{-1}\rho\br{-\mu p\fr{v}{F}\fr{\vt'}{\vt}+(p+1)\fr{\mu-1}{ n}}.}
Set
\eq{\ti{\rho}(t)=\max_{M}\rho(t,\cdot).}
Note that $\ti{\rho}$ is Lipschitz continuous, hence differentiable almost everywhere in $(0,\8)$ and at points of differentiability there holds
\eq{\dot{\ti{\rho}}(t)=\fr{\del \rho}{\del t}(t,x_{t}),}
where 
\eq{\rho(t,x_{t})=\ti{\rho}(t),}
cf. \cite[Lemma~6.3.2]{Gerhardt:/2006}. The original idea of this useful fact goes back to Hamilton \cite[Lemma~3.5]{Hamilton:/1986}. 
Choosing $\mu>0$ small enough, $\|\n u\|^{2}\vt'^{\mu}$ converges to zero due to \cref{ICF-grad-bound} and we obtain that for sufficiently large $t$, 
\eq{\dot{\ti{\rho}}(t)\leq 0}
 on the set $\{\ti{\rho}\geq 1\}$, provided $p<1$. Hence in this case $\rho$ is bounded. In case $p=1$ we set $\mu=0$ and obtain that for all $\e>0$ there exist $\de_{\e}>0$ and $T_{\e}$, such that for all $t\geq T_{\e}$ where $\ti{\rho}$ is differentiable, there holds
 \eq{\ti\rho(t)\geq \e\q\Ra\q \dot{\ti{\rho}}(t)< -\de_{\e}.}
 \cite[Lemma~4.2]{Scheuer:05/2015} implies $\limsup_{t\ra\8} \ti\rho\leq 0$. Hence
 \eq{\limsup_{t\ra \8}v\ka_{n}\fr{\vt}{\vt'}\leq 1}
 in both cases.
Now
 \eq{\sum_{i=1}^{n}\fr{\br{1-v\ka_{i}\fr{\vt}{\vt'}}}{nvF\fr{\vt}{\vt'}}=\fr{n-vH\fr{\vt}{\vt'}}{nvF\fr{\vt}{\vt'}}\leq \fr{n-vF\fr{\vt}{\vt'}}{nvF\fr{\vt}{\vt'}}\leq \fr{\vt'}{F\vt}-\fr{1}{n}\leq c\vt'^{-\mu}}
 and hence
 \eq{\label{ICF-LT-k-Decay-2}1-v\ka_{1}\fr{\vt}{\vt'}\leq c\vt'^{-\mu} +\sum_{i=2}^{n}\br{v\ka_{i}\fr{\vt}{\vt'}-1}.}
 The proof is complete. 
}

Now we are in the position to optimise the decay estimates. We start with the gradient.

\begin{lemma}\label{ICF-optGrad}
Under the assumption of \cref{ICF-Main}~(ii) the function
\eq{\label{ICF-optGrad-A}\tilde{z}=|\hat\n\p|^{2}\vt'^{2p}}
is uniformly bounded. If $p=1$,
then additionally there exist constants $c$ and $\al$ such that
\eq{|\hat\n\p|^{2}\leq ce^{-\al t}.}
\end{lemma}

\pf{If $\vt'$ is unbounded, using \cref{ICF-LT-k-Decay} we can rewrite the evolution of
\eq{z=f(u)|\hat\n\p|^{2}}
from \eqref{ICF-grad-bound-1} with 
\eq{f(u)=\vt'^{\g}(u)}
at maximal points as
\eq{\label{ICF-optGrad-1}\cL z&\leq\fr{2z}{F^{p+1}}\Big((p-1)\fr{\vt'}{\vt}vF-p\fr{\vt''}{\vt}F^{k}_{k}+\g\fr{p+1}{2}\fr{\vt''}{\vt'}\fr{F}{v}-\fr{\g np}{2}\fr{\vt''}{\vt}\Big)\\
&\hp{=}-\fr{2f}{\vt^{2}}\fr{p}{F^{p+1}}F^{k}_{l}\ti{g}^{lr}\hat{R}^{m}_{ikr}\p^{i}\p_{m}+\fr{p}{F^{p+1}}F^{kr}u_{k}u_{r}z\br{2\fr{f'^{2}}{f^{2}}-\fr{f''}{f}-\fr{f'}{f}\fr{\vt'}{\vt}}\\
&\leq \fr{2z}{F^{p+1}}\br{o(1)\fr{\vt'^{2}}{\vt^{2}}+n(p-1)\fr{\vt'^{2}}{\vt^{2}}-np\fr{\vt''}{\vt}+\fr{\g n}{2}\fr{\vt''}{\vt}}\\
&\hp{=}-\fr{2f}{\vt^{2}}\fr{p}{F^{p+1}}F^{k}_{l}\ti{g}^{lr}\hat{R}^{m}_{ikr}\p^{i}\p_{m}.
}
Since we want to bound $z$, it suffices to consider spatial maxima at which $z$ is positive. At such there holds
\eq{\label{ICF-optGrad-2}\cL z&\leq\fr{2z}{F^{p+1}}\fr{\vt'^{2}}{\vt^{2}}\Big(o(1)+n(p-1)-n\br{p-\fr{\g}{2}}\fr{\vt''\vt}{\vt'^{2}}\\
					&\hp{\fr{2z}{F^{p+1}}\fr{\vt'^{2}}{\vt^{2}}\Big(o(1)}-\fr{p}{\vt'^{2}}\widehat{\Rc}\br{\fr{\hat\n\p}{|\hat\n\p|},\fr{\hat\n\p}{|\hat\n\p|}}\Big).
}
In case $p<1$ with $\g=2p$, the right hand side is eventually negative for large $t$, since only the case of unbounded $\vt'$ has to be considered to prove the first statement. In case $p=1$ we put $\g=0$ and use the first estimate in \eqref{ICF-optGrad-1} if $\vt'$ is bounded, whereas if $\vt'$ is unbounded we use \eqref{ICF-optGrad-2}, to get
\eq{\cL z\leq -\de z}
for some $\de$ and large times. The exponential decay follows. To prove the remaining claim, we evaluate \eqref{ICF-grad-bound-2} with 
\eq{f=\vt'^{2},\q p=1}
 and see
\eq{\cL z\leq -\fr{2}{F^{2}}\fr{\vt''}{\vt}\br{1+1-ce^{-\al t}-2}F^{k}_{k}z+\fr{c}{F^{2}}\fr{\vt'^{2}}{\vt^{2}}e^{-\al t}z.}
Hence the function
\eq{\bar{z}(t)=\max_{\cS_{0}}z(t,\cdot)}
satisfies
\eq{\dot{\bar{z}}\leq ce^{-\al t}\bar{z}}
and is thus bounded.
}

We optimise the convergence rate of the rescaled principal curvatures.

\begin{lemma}\label{ICF-opt-k-Decay}
Under the assumptions of \cref{ICF-Main}~(ii), along \eqref{ICF} there exists a constant $c$, such that for all $1\leq i\leq n$, the $i$-th rescaled principal curvature satisfies
\eq{\left|v\kappa_{i}\fr{\vt}{\vt'}-1\right|\leq \fr{ct}{\vt'^{p(p+1)}},}
where we may drop the $t$-factor if $p<1$ or if $\vt'$ is bounded.
\end{lemma}

\pf{Only the case that $\vt'$ is unbounded has to be considered. 

(i)~First we optimise the decay in \cref{ICF-LT-1/F-Decay}.
Using the optimal gradient estimates \cref{ICF-optGrad}, we see from \eqref{ICF-LT-1/F-Decay-1} that \cref{ICF-LT-1/F-Decay} holds with any $\mu<p(p+1)$, if $c$ is allowed to depend on a lower bound of $p(p+1)-\mu$.

Now consider the function $\rho$ defined in the proof of \cref{ICF-LT-k-Decay} and obtain from \eqref{ICF-LT-k-Decay-1} with $\mu<p(p+1)$ at points where $\rho\geq 1$ that
\eq{&\dot{\rho}-\fr{p}{F^{p+1}}F^{ij}\rho_{;ij}\\
	\leq~&\fr{\vt'}{\vt}\fr{p+1}{vF^{p}}\br{-\fr{1-p}{p+1}\rho+c\|\n u\|^{2}\vt'^{\mu}+\fr{c\vt'^{\mu}}{\vt'^{2}}\|\n u\|}\\
	+&~\fr{n}{F^{p}}\fr{\vt''}{\vt'}v^{-1}\rho\br{\fr{-p\mu}{n}+o(1)+(p+1)\fr{\mu-1}{ n}}\\
	=~&\fr{\vt'}{\vt}\fr{p+1}{vF^{p}}\br{-\fr{1-p}{p+1}\rho+c\|\n u\|^{2}\vt'^{\mu}+\fr{c\vt'^{\mu}}{\vt'^{2}}\|\n u\|}\\
	+&~\fr{1}{F^{p}}\fr{\vt''}{\vt'}v^{-1}\rho\br{\mu-(p+1)+o(1)}\\
	<~&0
	}
in case $p<1$ for large times. In case $p=1$ the right hand side of this inequality eventually decays exponentially and thus
\eq{\rho\leq c}
in both cases.
Hence for any $\mu<p(p+1)$ we have
\eq{v\ka_{n}\fr{\vt}{\vt'}-1\leq \fr{c_{\mu}}{\vt'^{\mu}}.}

 Now putting $\mu=p(p+1)$ in \eqref{ICF-LT-1/F-Decay-1} we see that the function $\rho$ defined in the proof of \cref{ICF-LT-1/F-Decay} satisfies at positive maximal points with $\rho\geq 1$
\eq{&\dot{\rho}-\fr{p}{F^{p+1}}F^{ij}\rho_{;ij}\\
	\leq~&\fr{np(p+1)}{F^{p+1}}\fr{\vt''}{\vt}e^{\fr{p-1}{p}z}\Big(-\rho+(p+1)\fr{F}{n}\fr{\vt}{\vt'}v^{-1}e^{\fr{1-p}{p}z}\rho-p\rho+c\|\n u\|^{2}\vt'^{p(p+1)}\Big)\\
	+&~\fr{np(1-p)}{F^{p+1}}\fr{\vt'^{2}}{\vt^{2}}e^{\fr{p-1}{p}z}\br{-\rho+c\|\n u\|^{2}\vt'^{p(p+1)}}\\
	<~&\fr{np(p+1)}{F^{p+1}}\fr{\vt''}{\vt}e^{\fr{p-1}{p}z}\Big(c\vt'^{-p}\rho+c\|\n u\|^{2}\vt'^{p(p+1)}\Big)\\
	+&~\fr{np(1-p)}{F^{p+1}}\fr{\vt'^{2}}{\vt^{2}}e^{\fr{p-1}{p}z}\br{-\rho+c\|\n u\|^{2}\vt'^{p(p+1)}}.
}
In case $p<1$ we use
\eq{\fr{\vt''}{\vt}\leq c\fr{\vt'^{2}}{\vt^{2}}}
to absorb every decaying term into $-\rho$ in the second line.
In case $p=1$ we use 
\eq{\vt'^{-p}\rho\leq c}
to conclude
\eq{\dot{\rho}-\fr{p}{F^{p+1}}F^{ij}\rho_{;ij}\leq c.}
Hence we obtain
\eq{\fr{\vt'}{\vt}{\fr{1}{F}}-\fr{1}{n}\leq \fr{ct}{\vt'^{p(p+1)}},}
and the same without the $t$-factor in case $p<1$.

(ii)~In the second step we optimise the convergence rate in \cref{ICF-LT-k-Decay}. Therefore we consider the function $\rho$ defined in that proof and obtain from \eqref{ICF-LT-k-Decay-1} with $\mu=p(p+1)$ at points where $\rho\geq 1$ that
\eq{&\dot{\rho}-\fr{p}{F^{p+1}}F^{ij}\rho_{;ij}\\
	\leq~&\fr{\vt'}{\vt}\fr{p+1}{vF^{p}}\br{-\fr{1-p}{p+1}\rho+c\|\n u\|^{2}\vt'^{p(p+1)}+\fr{c\vt'^{p(p+1)}}{\vt'^{2}}\|\n u\|}\\
	+&~\fr{n}{F^{p}}\fr{\vt''}{\vt'}v^{-1}\rho\br{\fr{-p^{2}(p+1)}{n}+o(1)+(p+1)\fr{p(p+1)-1}{ n}}\\
	=~&\fr{\vt'}{\vt}\fr{p+1}{vF^{p}}\br{-\fr{1-p}{p+1}\rho+c\|\n u\|^{2}\vt'^{p(p+1)}+\fr{c\vt'^{p(p+1)}}{\vt'^{2}}\|\n u\|}\\
	+&~\fr{1}{F^{p}}\fr{\vt''}{\vt'}v^{-1}\rho\br{p^{2}-1+o(1)}\\
	\leq~&0
	}
in case $p<1$ for large times. In case $p=1$ the right hand side of this inequality is bounded and thus
\eq{\rho\leq ct}
in this case.
Estimating \eqref{ICF-LT-k-Decay-2} with the optimised bounds completes the proof.
}

We finish the proof of \cref{ICF-Main} by proving the final statement about the exponential decay in item (ii).
The function 
\eq{z=\log\vt(u)-\fr{t}{n}}
defined on $[0,\8)\x \cS_{0}$ satisfies
\eq{\dot{z}=\fr{v\vt'}{\vt F(\cW)}-\fr{1}{n}=\fr{v}{F\br{\fr{1}{v}\de^{i}_{j}+\fr{1}{v^{3}\vt^{2}}u^{i}u_{j}-\fr{1}{v\vt'\vt}\ti{g}^{ik}u_{kj}}}-\fr 1n=G(y,z,\hat\n z,\hat\n^{2}z).}
Hence
\eq{\fr{\del G}{\del z_{ij}}=\fr{\vt'^{-2}}{F^{2}\br{\fr{\vt}{\vt'}\cW}}F^{i}_{k}\ti{g}^{kj},}
which is uniformly elliptic, since $\vt'$ is globally bounded. $z$ is uniformly bounded, as can be seen similarly as in \cite[Prop.~3.1, Lemma~3.2]{Scheuer:01/2017}. Furthermore
\eq{|\hat\n z|\leq c|\hat\n\p|\leq ce^{-\al t}}
and 
\eq{|\hat\n^{2}\p|\leq c.}
Applying the regularity results of Krylov and Safonov as well as Schauder theory, we obtain uniform $C^{m}$-bounds for $z$. Due to interpolation we get
\eq{|\hat\n^{2}\p|\leq ce^{-\al t},}
which implies \eqref{ICF-Umb} with $t$ replaced by $e^{-\al t}$.

\begin{rem}\label{HigherOrder}
The previous argument is precisely the way to deduce a uniform bound on the rescaled principal curvatures for the inverse mean curvature flow, when $\vt'$ is bounded, as it was performed in \cite{Mullins:10/2016, Zhou:06/2017}.  
The crucial point is here, that one does not need curvature estimates to have $F^{ij}$ uniformly elliptic. One only needs a bound on the rescaled speed
\eq{\ti{H}=\fr{\vt}{\vt'}H.}
Then the above argumentation applies.

\end{rem}

\bibliographystyle{amsplain}
\bibliography{/Users/J_Mac/Documents/Uni/TexTemplates/bibliography}

\end{document}